\documentclass{article}

\usepackage{latexsym}
\usepackage{amsmath}
\usepackage{amssymb}
\usepackage{amsthm}

\usepackage{moreverb}

\newtheorem{alg}{Algorithm}

\usepackage{parskip}

\newcommand{\jjbls}[0]{\small\baselineskip0.75em}
\newcommand{\jjterm}[1]{\emph{#1}}
\newcommand{\eqq}[1]{``{#1}''}
\newcommand{\codesize}[0]{\small\baselineskip0.65em}
\newcommand{\fileref}[1]{\texttt{#1}}

\newcommand{\jjseqref}[1]{\href{http://oeis.org/#1}{#1}}

\usepackage[svgnames]{xcolor}

\usepackage{ifpdf}
\ifpdf
\usepackage[hyperindex,pdftex,backref=page,breaklinks=true,colorlinks=true,linkcolor=DarkBlue,citecolor=DarkBlue,urlcolor=DarkBlue]{hyperref}
\hypersetup{%
pdfauthor={Joerg Arndt},%
pdftitle={Subset-lex: did we miss an order?},%
pdfsubject={combinatorial generation},%
pdfkeywords={combinatorial generation, subset-lex order, lexicographic order, Gray code, loopless generation, O(1) algorithms},%
bookmarksopen,
bookmarksopenlevel={0},
}
\else
\usepackage[hypertex,hyperindex,backref=page,breaklinks=true,colorlinks=true,linkcolor=DarkBlue,citecolor=DarkBlue,urlcolor=DarkBlue]{hyperref}%
\fi

\begin{document}
\bibliographystyle{plain}
\title{Subset-lex: did we miss an order?}
\author{J\"{o}rg Arndt}
\date{January 12, 2015}
\pagestyle{plain} \maketitle

\begin{abstract}\noindent

\end{abstract}
We generalize a well-known algorithm for the generation of
all subsets of a set in lexicographic order with respect to
the sets as lists of elements (subset-lex order).
We obtain algorithms for various combinatorial objects
such as the subsets of a multiset,
compositions and partitions represented as lists of parts,
and for certain restricted growth strings.
The algorithms are often loopless and require at most
one extra variable for the computation of the next object.
The performance of the algorithms is very
competitive even when not loopless.
A Gray code corresponding to the subset-lex order
and a Gray code for compositions that was
found during this work are described.

Appendix II about ranking and unranking methods for mixed radix words
added January 2, 2024 at the very end. The text is otherwise unchanged.

{
\baselineskip3.0mm
\tableofcontents
}

\section{Two lexicographic orders for binary words}

%
\begin{figure}
{\jjbls
\begin{verbatim}
 0:  [ . . . . . ]  {  }                 [ . . . . . ]  {  }
 1:  [ 1 . . . . ]  { 0 }                [ . . . . 1 ]  { 4 }
 2:  [ 1 1 . . . ]  { 0, 1 }             [ . . . 1 . ]  { 3 }
 3:  [ 1 1 1 . . ]  { 0, 1, 2 }          [ . . . 1 1 ]  { 3, 4 }
 4:  [ 1 1 1 1 . ]  { 0, 1, 2, 3 }       [ . . 1 . . ]  { 2 }
 5:  [ 1 1 1 1 1 ]  { 0, 1, 2, 3, 4 }    [ . . 1 . 1 ]  { 2, 4 }
 6:  [ 1 1 1 . 1 ]  { 0, 1, 2, 4 }       [ . . 1 1 . ]  { 2, 3 }
 7:  [ 1 1 . 1 . ]  { 0, 1, 3 }          [ . . 1 1 1 ]  { 2, 3, 4 }
 8:  [ 1 1 . 1 1 ]  { 0, 1, 3, 4 }       [ . 1 . . . ]  { 1 }
 9:  [ 1 1 . . 1 ]  { 0, 1, 4 }          [ . 1 . . 1 ]  { 1, 4 }
10:  [ 1 . 1 . . ]  { 0, 2 }             [ . 1 . 1 . ]  { 1, 3 }
11:  [ 1 . 1 1 . ]  { 0, 2, 3 }          [ . 1 . 1 1 ]  { 1, 3, 4 }
12:  [ 1 . 1 1 1 ]  { 0, 2, 3, 4 }       [ . 1 1 . . ]  { 1, 2 }
13:  [ 1 . 1 . 1 ]  { 0, 2, 4 }          [ . 1 1 . 1 ]  { 1, 2, 4 }
14:  [ 1 . . 1 . ]  { 0, 3 }             [ . 1 1 1 . ]  { 1, 2, 3 }
15:  [ 1 . . 1 1 ]  { 0, 3, 4 }          [ . 1 1 1 1 ]  { 1, 2, 3, 4 }
16:  [ 1 . . . 1 ]  { 0, 4 }             [ 1 . . . . ]  { 0 }
17:  [ . 1 . . . ]  { 1 }                [ 1 . . . 1 ]  { 0, 4 }
18:  [ . 1 1 . . ]  { 1, 2 }             [ 1 . . 1 . ]  { 0, 3 }
19:  [ . 1 1 1 . ]  { 1, 2, 3 }          [ 1 . . 1 1 ]  { 0, 3, 4 }
20:  [ . 1 1 1 1 ]  { 1, 2, 3, 4 }       [ 1 . 1 . . ]  { 0, 2 }
21:  [ . 1 1 . 1 ]  { 1, 2, 4 }          [ 1 . 1 . 1 ]  { 0, 2, 4 }
22:  [ . 1 . 1 . ]  { 1, 3 }             [ 1 . 1 1 . ]  { 0, 2, 3 }
23:  [ . 1 . 1 1 ]  { 1, 3, 4 }          [ 1 . 1 1 1 ]  { 0, 2, 3, 4 }
24:  [ . 1 . . 1 ]  { 1, 4 }             [ 1 1 . . . ]  { 0, 1 }
25:  [ . . 1 . . ]  { 2 }                [ 1 1 . . 1 ]  { 0, 1, 4 }
26:  [ . . 1 1 . ]  { 2, 3 }             [ 1 1 . 1 . ]  { 0, 1, 3 }
27:  [ . . 1 1 1 ]  { 2, 3, 4 }          [ 1 1 . 1 1 ]  { 0, 1, 3, 4 }
28:  [ . . 1 . 1 ]  { 2, 4 }             [ 1 1 1 . . ]  { 0, 1, 2 }
29:  [ . . . 1 . ]  { 3 }                [ 1 1 1 . 1 ]  { 0, 1, 2, 4 }
30:  [ . . . 1 1 ]  { 3, 4 }             [ 1 1 1 1 . ]  { 0, 1, 2, 3 }
31:  [ . . . . 1 ]  { 4 }                [ 1 1 1 1 1 ]  { 0, 1, 2, 3, 4 }
\end{verbatim}
}
\caption{\label{fig:subsetlex}
Two lexicographic orders for the subsets of a 5-element set:
ordering using the sets as lists of elements (left)
and ordering using the characteristic words (right).
Dots are used to denote zeros in the characteristic words.}
\end{figure}
%

Two simple ways to represent subsets of a set $\{0,1,2,\ldots,n-1\}$ are
as lists of elements and as the characteristic word, the binary word with
ones at the positions corresponding to the elements in the subset.
Sorting all subsets by list of characteristic words gives the right columns of
figure \ref{fig:subsetlex}, sorting by the list of elements gives the left
columns.
We will concern ourselves with the generalization of the ordering by the lists
which we call \jjterm{subset-lex} order.

The algorithms for computing the successor or predecessor
of a (nonempty) subset $S$ in subset-lex order
are well-known\footnote{For example,
 an algorithm for the $k$-subsets of an $n$-set is given
 in \cite[Algorithm LEXSUB, p.18]{nw-combalg},
 see section \ref{sect:ksubset} for an implementation.}.
In the following let $z$ be the last, and $y$ the next to last element in $S$.
We call an element minimal or maximal if it is respectively
the smallest or greatest element in the superset
%
\begin{alg}[Next-SL2]\label{alg:next-sl2}
Compute the successor of the subset $S$.
\end{alg}
\begin{enumerate}
\item If there is just one element, and it is maximal, stop.
\item If $z$ is not maximal, append $z+1$.
\item Otherwise remove $z$ and $y$, then append $y+1$.
\end{enumerate}
\begin{alg}[Prev-SL2]
Compute the predecessor of the subset $S$.
\end{alg}
\begin{enumerate}
\item If there is just one element, and it is minimal, stop.
\item If $z-1 \in S$, remove $z$.
\item Otherwise remove $z$, append $z-1$ and the maximal element.
\end{enumerate}
%

C++ implementations of all algorithms discussed here are given in the FXT
library \cite{fxt}.
Here and elsewhere we give slightly simplified versions of the actual code,
like omitting modifiers such as \texttt{public}, \texttt{protected}, and
\texttt{private}.  The type \texttt{ulong} is short for \texttt{unsigned long}.
{\codesize
\begin{listing}{1}
// FILE:  src/comb/subset-lex.h
class subset_lex
// Nonempty subsets of the set {0,1,2,...,n-1} in subset-lex order.
// Representation as list of parts.
// Loopless generation.
{
    ulong n_;   // number of elements in set, should have n>=1
    ulong k_;   // index of last element in subset
    ulong n1_;  // == n - 1 for n >=1, and == 0 for n==0
    ulong *x_;  // x[0...k-1]:  subset of {0,1,2,...,n-1}

    subset_lex(ulong n)
    // Should have n>=1,
    // for n==0 one set with the element zero is generated.
    {
        n_ = n;
        n1_ = (n_ ? n_ - 1 : 0 );
        x_ = new ulong[n_ + (n_==0)];
        first();
    }
\end{listing}
}%
We will always use the method names
\texttt{first()} and \texttt{last()}
for the first and last element in
the respective list of combinatorial objects.
{\codesize
\begin{listing}{1}
    ulong first()
    {
        k_ = 0;
        x_[0] = 0;
        return  k_ + 1;
    }

    ulong last()
    {
        k_ = 0;
        x_[0] = n1_;
        return  k_ + 1;
    }
\end{listing}
}%

The methods for computing the successor and predecessor
will always be called \texttt{next()} and \texttt{prev()}
respectively.
{\codesize
\begin{listing}{1}
    ulong next()
    // Return number of elements in subset.
    // Return zero if current is last.
    {
        if ( x_[k_] == n1_ )  // last element is max?
        {
            if ( k_==0 )  return 0;

            --k_;      // remove last element
            x_[k_] += 1;  // increment last element
        }
        else  // add next element from set:
        {
            ++k_;
            x_[k_] = x_[k_-1] + 1;
        }

        return  k_ + 1;
    }

    ulong prev()
    // Return number of elements in subset.
    // Return zero if current is first.
    {
        if ( k_ == 0 )  // only one element ?
        {
            if ( x_[0]==0 )  return 0;

            x_[0] -= 1;  // decrement last (and only) element
            ++k_;
            x_[k_] = n1_;  // append maximal element
        }
        else
        {
            if ( x_[k_] == x_[k_-1] + 1 )  --k_;  // remove last element
            else
            {
                x_[k_] -= 1;  // decrement last element
                ++k_;
                x_[k_] = n1_;  // append maximal element
            }
        }

        return  k_ + 1;
    }
\end{listing}
}%
A program demonstrating usage of the class is
{\codesize
\begin{listing}{1}
// FILE:  demo/comb/subset-lex-demo.cc
    ulong n = 6;
    subset_lex S(n);
    do
    {
        // visit subset
    }
    while ( S.next() );
\end{listing}
}%
%
%

\section{Subset-lex order for multisets}

%
\begin{figure}
\vspace*{-5mm}
{\jjbls
\begin{verbatim}
       0:    [ . . . . ]    {  }
       1:    [ 1 . . . ]    { 0 }
       2:    [ 1 1 . . ]    { 0, 1 }
       3:    [ 1 2 . . ]    { 0, 1, 1 }
       4:    [ 1 2 1 . ]    { 0, 1, 1, 2 }
       5:    [ 1 2 2 . ]    { 0, 1, 1, 2, 2 }
       6:    [ 1 2 2 1 ]    { 0, 1, 1, 2, 2, 3 }
       7:    [ 1 2 2 2 ]    { 0, 1, 1, 2, 2, 3, 3 }
       8:    [ 1 2 2 3 ]    { 0, 1, 1, 2, 2, 3, 3, 3 }
       9:    [ 1 2 1 1 ]    { 0, 1, 1, 2, 3 }
      10:    [ 1 2 1 2 ]    { 0, 1, 1, 2, 3, 3 }
      11:    [ 1 2 1 3 ]    { 0, 1, 1, 2, 3, 3, 3 }
      12:    [ 1 2 . 1 ]    { 0, 1, 1, 3 }
      13:    [ 1 2 . 2 ]    { 0, 1, 1, 3, 3 }
      14:    [ 1 2 . 3 ]    { 0, 1, 1, 3, 3, 3 }
      15:    [ 1 1 1 . ]    { 0, 1, 2 }
      16:    [ 1 1 2 . ]    { 0, 1, 2, 2 }
      17:    [ 1 1 2 1 ]    { 0, 1, 2, 2, 3 }
      18:    [ 1 1 2 2 ]    { 0, 1, 2, 2, 3, 3 }
      19:    [ 1 1 2 3 ]    { 0, 1, 2, 2, 3, 3, 3 }
      20:    [ 1 1 1 1 ]    { 0, 1, 2, 3 }
      21:    [ 1 1 1 2 ]    { 0, 1, 2, 3, 3 }
      22:    [ 1 1 1 3 ]    { 0, 1, 2, 3, 3, 3 }
      23:    [ 1 1 . 1 ]    { 0, 1, 3 }
      24:    [ 1 1 . 2 ]    { 0, 1, 3, 3 }
      25:    [ 1 1 . 3 ]    { 0, 1, 3, 3, 3 }
      26:    [ 1 . 1 . ]    { 0, 2 }
      27:    [ 1 . 2 . ]    { 0, 2, 2 }
      28:    [ 1 . 2 1 ]    { 0, 2, 2, 3 }
      29:    [ 1 . 2 2 ]    { 0, 2, 2, 3, 3 }
      30:    [ 1 . 2 3 ]    { 0, 2, 2, 3, 3, 3 }
      31:    [ 1 . 1 1 ]    { 0, 2, 3 }
      32:    [ 1 . 1 2 ]    { 0, 2, 3, 3 }
      33:    [ 1 . 1 3 ]    { 0, 2, 3, 3, 3 }
      34:    [ 1 . . 1 ]    { 0, 3 }
      35:    [ 1 . . 2 ]    { 0, 3, 3 }
      36:    [ 1 . . 3 ]    { 0, 3, 3, 3 }
      37:    [ . 1 . . ]    { 1 }
      38:    [ . 2 . . ]    { 1, 1 }
      39:    [ . 2 1 . ]    { 1, 1, 2 }
      40:    [ . 2 2 . ]    { 1, 1, 2, 2 }
      41:    [ . 2 2 1 ]    { 1, 1, 2, 2, 3 }
      42:    [ . 2 2 2 ]    { 1, 1, 2, 2, 3, 3 }
      43:    [ . 2 2 3 ]    { 1, 1, 2, 2, 3, 3, 3 }
      44:    [ . 2 1 1 ]    { 1, 1, 2, 3 }
      45:    [ . 2 1 2 ]    { 1, 1, 2, 3, 3 }
      46:    [ . 2 1 3 ]    { 1, 1, 2, 3, 3, 3 }
      47:    [ . 2 . 1 ]    { 1, 1, 3 }
      48:    [ . 2 . 2 ]    { 1, 1, 3, 3 }
      49:    [ . 2 . 3 ]    { 1, 1, 3, 3, 3 }
      50:    [ . 1 1 . ]    { 1, 2 }
      51:    [ . 1 2 . ]    { 1, 2, 2 }
      52:    [ . 1 2 1 ]    { 1, 2, 2, 3 }
      53:    [ . 1 2 2 ]    { 1, 2, 2, 3, 3 }
      54:    [ . 1 2 3 ]    { 1, 2, 2, 3, 3, 3 }
      55:    [ . 1 1 1 ]    { 1, 2, 3 }
      56:    [ . 1 1 2 ]    { 1, 2, 3, 3 }
      57:    [ . 1 1 3 ]    { 1, 2, 3, 3, 3 }
      58:    [ . 1 . 1 ]    { 1, 3 }
      59:    [ . 1 . 2 ]    { 1, 3, 3 }
      60:    [ . 1 . 3 ]    { 1, 3, 3, 3 }
      61:    [ . . 1 . ]    { 2 }
      62:    [ . . 2 . ]    { 2, 2 }
      63:    [ . . 2 1 ]    { 2, 2, 3 }
      64:    [ . . 2 2 ]    { 2, 2, 3, 3 }
      65:    [ . . 2 3 ]    { 2, 2, 3, 3, 3 }
      66:    [ . . 1 1 ]    { 2, 3 }
      67:    [ . . 1 2 ]    { 2, 3, 3 }
      68:    [ . . 1 3 ]    { 2, 3, 3, 3 }
      69:    [ . . . 1 ]    { 3 }
      70:    [ . . . 2 ]    { 3, 3 }
      71:    [ . . . 3 ]    { 3, 3, 3 }
\end{verbatim}
}
\caption{\label{fig:msubsetlex}
Subsets of the multiset $\{0^1, 1^2, 2^2, 3^3 \}$ in subset-lex order.
Dots are used to denote zeros in the (generalized) characteristic words.
}
\end{figure}
%

For the internal representation of the subsets of a multiset, we could use lists
of pairs $(e, m)$ where $e$ is the element and $m$ its multiplicity.  This
choice would lead to loopless algorithms as will become clear in a moment.  The
disadvantage of this representation, however, is that the generalizations we
will find would have a slightly more complicated form.
We will instead use generalized characteristic words,
where nonzero entries can be at most the multiplicity
of the element in question.
For a multiset $\{0^{m(0)},\, 1^{m(1)},\, \ldots,\, k^{m(k)}\}$
these are the mixed radix numbers with radix vector
$[m(0)+1,\, m(1)+1,\, \ldots,\, m(k)+1]$.
The subsets of the set $\{0^1,\, 1^2,\, 2^2,\, 3^3 \}$ in subset-lex order
are shown in figure \ref{fig:msubsetlex}.

The algorithms for computing the successor or predecessor
of a (nonempty) subset $S$ of a multiset in subset-lex are now given.
In the following let $z$ be the last element in $S$.
%
\begin{samepage}
\begin{alg}[Next-SL]\label{alg:next-sl}
Compute the successor of the subset $S$.
\end{alg}
\begin{enumerate}
\item If the multiplicity of $z$ is not maximal, increase it;  return.
\item If $z$ is not maximal, append $z+1$ with multiplicity 1;  return.
\item If $z$ is the only element, stop.
\item Remove $z$ and decrease the multiplicity of next to last nonzero element $y$,
  then append $y+1$ with multiplicity 1.
\end{enumerate}
\end{samepage}
The description omits the case of the empty set,
the implementation takes care of this case by intially pointing
to the leftmost digit of the characteristic word, which is zero.
With the first call to the routine it is changed to 1.
\begin{alg}[Prev-SL]\label{alg:prev-sl}
Compute the predecessor of the subset $S$.
\end{alg}
\begin{enumerate}
\item If the set is empty, stop.
\item Decrease the multiplicity of $z$.
\item If the new multiplicity is nonzero, return.
\item Increase the multiplicity of $z-1$ and append the maximal element.
\end{enumerate}

Clearly both algorithms are loopless for the representation using a list of pairs.
With the representation as characteristic words algorithm \ref{alg:next-sl}
involves a loop in the last step, a scan for the next to last nonzero element.
This could be prevented by maintaining an additional list of positions where the
characteristic word is nonzero.
Algorithm \ref{alg:prev-sl} does stay loopless, the position of
the maximal element does not need to be sought.
Our implementations will always use a variable
holding the position of the last nonzero position
in the characteristic word.
This is called the \eqq{current track} in the implementations.

We will not make the equivalents of algorithm \ref{alg:next-sl} loopless, as
the maintenance of the additional list would render algorithm \ref{alg:prev-sl}
slower (at least as long as mixed calls to both shall be allowed),
but see section \ref{sect:loopless-next} for one such example.

A step in either direction has either one or three positions in the
characteristic word changed.  The number of transitions with one change
are more frequent with higher multiplicities.  The worst case occurs
if all multiplicities are one (corresponding to binary words), when
(about) half of the steps involve only one transition.

The C++ implementation uses the sentinel technique
to reduce the number of conditional branches.
{\codesize
\begin{listing}{1}
// FILE:  src/comb/mixedradix-subset-lex.h
class mixedradix_subset_lex
{
    ulong n_;    // number of digits (n kinds of elements in multiset)
    ulong tr_;   // aux: current track
    ulong *a_;   // digits of mixed radix number
                 //  (multiplicities in subset).
    ulong *m1_;  // nines (radix minus one) for each digit
                 //  (multiplicity of kind k in superset).

    mixedradix_subset_lex(ulong n, ulong mm, const ulong *m=0)
    {
        n_ = n;
        a_ = new ulong[n_+2];  // two sentinels, one left, one right
        a_[0] = 1;  a_[n_+1] = 1;
        ++a_;  // nota bene
        m1_ = new ulong[n_+2];
        m1_[0] = 0;  m1_[n_+1] = 0;  // sentinel with n==0
        ++m1_;  // nota bene

        mixedradix_init(n_, mm, m, m1_);  // set up m1_[], omitted

        first();
    }
\end{listing}
}%
The omitted routine \texttt{mixedradix\_init()}
sets all elements of the array \texttt{m1[]} to \texttt{mm}
(fixed radix case),
unless the pointer \texttt{m} is nonzero, then the
elements behind \texttt{m} are copied into \texttt{m1[]}
(mixed radix case).

{\codesize
\begin{listing}{1}
    void first()
    {
        for (ulong k=0; k<n_; ++k)  a_[k] = 0;
        tr_ = 0;  // start by looking at leftmost (zero) digit
    }

    void last()
    {
        for (ulong k=0; k<n_; ++k)  a_[k] = 0;
        ulong n1 = ( n_ ? n_ - 1 : 0 );
        a_[n1] = m1_[n1];
        tr_ = n1;
    }
\end{listing}
}%

Note the \texttt{while}-loop in the \texttt{next()} method.
{\codesize
\begin{listing}{1}
    bool next()
    {
        ulong j = tr_;
        if ( a_[j] < m1_[j] )  // easy case 1: increment
        {
            a_[j] += 1;
            return true;
        }

        // here a_[j] == m1_[j]
        if ( j+1 < n_ )  // easy case 2: append (move track to the right)
        {
            ++j;
            a_[j] = 1;
            tr_ = j;
            return true;
        }

        a_[j] = 0;

        // find first nonzero digit to the left:
        --j;
        while ( a_[j] == 0 )  { --j; }  // may read sentinel a_[-1]

        if ( (long)j < 0 )  return false;  // current is last

        a_[j] -= 1;  // decrement digit to the left
        ++j;
        a_[j] = 1;
        tr_ = j;
        return true;
    }
\end{listing}
}%
The method \texttt{prev()} is indeed loopless:
{\codesize
\begin{listing}{1}
    bool prev()
    {
        ulong j = tr_;
        if ( a_[j] > 1 )  // easy case 1: decrement
        {
            a_[j] -= 1;
            return true;
        }
        else
        {
            if ( tr_ == 0 )
            {
                if ( a_[0] == 0 )  return false;  // current is first
                a_[0] = 0;  // now word is first (all zero)
                return true;
            }

            a_[j] = 0;

            --j;  // now looking at next track to the left
            if ( a_[j] == m1_[j] )  // easy case 2: move track to left
            {
                tr_ = j;  // move track one left
            }
            else
            {
                a_[j] += 1;  // increment digit to the left
                j = n_ - 1;
                a_[j] = m1_[j];  // set rightmost digit = nine
                tr_ = j;  // move to rightmost track
            }
            return true;
        }
    }
\end{listing}
}%
%

%

The performance of the generator is quite satisfactory.
A system based on an AMD \texttt{Phenom(tm) II X4 945} processor
clocked at 3.0 GHz and the GCC C++ compiler version 4.9.0 \cite{GCC}
were used for measuring.
The updates using \texttt{next()} cost about
 12 cycles for base 2 (the worst case) and
 9 cycles for base 16.
Using \texttt{prev()} takes about
 8.5 cycles with base 2 and
 4.2 cycles with base 16,
the latter figure corresponding to a rate of 780 million subsets per second.

These and all following such figures are average values, obtained by measuring
the total time for the generation of all words of a certain size and dividing
by the number of words.


\section{Loopless computation of the successor for non-adjacent forms}\label{sect:loopless-next}

%
\begin{figure}
{\jjbls
\begin{verbatim}
           colex         Gray code       subset-lex     iset
   0:    [ . . . . ]    [ 4 . 2 . ]      [ . . . . ]    {  }
   1:    [ 1 . . . ]    [ 3 . 2 . ]      [ 1 . . . ]    { 0 }
   2:    [ 2 . . . ]    [ 2 . 2 . ]      [ 2 . . . ]    { 0 }
   3:    [ 3 . . . ]    [ 1 . 2 . ]      [ 3 . . . ]    { 0 }
   4:    [ 4 . . . ]    [ . . 2 . ]      [ 4 . . . ]    { 0 }
   5:    [ . 1 . . ]    [ . . 1 . ]      [ 4 . 1 . ]    { 0, 2 }
   6:    [ . 2 . . ]    [ 1 . 1 . ]      [ 4 . 2 . ]    { 0, 2 }
   7:    [ . 3 . . ]    [ 2 . 1 . ]      [ 4 . . 1 ]    { 0, 3 }
   8:    [ . . 1 . ]    [ 3 . 1 . ]      [ 3 . 1 . ]    { 0, 2 }
   9:    [ 1 . 1 . ]    [ 4 . 1 . ]      [ 3 . 2 . ]    { 0, 2 }
  10:    [ 2 . 1 . ]    [ 4 . . . ]      [ 3 . . 1 ]    { 0, 3 }
  11:    [ 3 . 1 . ]    [ 3 . . . ]      [ 2 . 1 . ]    { 0, 2 }
  12:    [ 4 . 1 . ]    [ 2 . . . ]      [ 2 . 2 . ]    { 0, 2 }
  13:    [ . . 2 . ]    [ 1 . . . ]      [ 2 . . 1 ]    { 0, 3 }
  14:    [ 1 . 2 . ]    [ . . . . ]      [ 1 . 1 . ]    { 0, 2 }
  15:    [ 2 . 2 . ]    [ . 1 . . ]      [ 1 . 2 . ]    { 0, 2 }
  16:    [ 3 . 2 . ]    [ . 2 . . ]      [ 1 . . 1 ]    { 0, 3 }
  17:    [ 4 . 2 . ]    [ . 3 . . ]      [ . 1 . . ]    { 1 }
  18:    [ . . . 1 ]    [ . 3 . 1 ]      [ . 2 . . ]    { 1 }
  19:    [ 1 . . 1 ]    [ . 2 . 1 ]      [ . 3 . . ]    { 1 }
  20:    [ 2 . . 1 ]    [ . 1 . 1 ]      [ . 3 . 1 ]    { 1, 3 }
  21:    [ 3 . . 1 ]    [ . . . 1 ]      [ . 2 . 1 ]    { 1, 3 }
  22:    [ 4 . . 1 ]    [ 1 . . 1 ]      [ . 1 . 1 ]    { 1, 3 }
  23:    [ . 1 . 1 ]    [ 2 . . 1 ]      [ . . 1 . ]    { 2 }
  24:    [ . 2 . 1 ]    [ 3 . . 1 ]      [ . . 2 . ]    { 2 }
  25:    [ . 3 . 1 ]    [ 4 . . 1 ]      [ . . . 1 ]    { 3 }
\end{verbatim}
}
\caption{\label{fig:naf}
Mixed radix numbers of length 4 in falling factorial base that are non-adjacent forms (NAF),
in co-lexicographic order, minimal-change order (Gray code),
and subset-lex order (with the sets of positions of nonzero digits).
}
\end{figure}
%

%
\begin{figure}
{\jjbls
\begin{verbatim}
   0:    [ 4 . 2 . ]    (0,  4) (1,  3) (2)       [ 4 3 2 1 0 ]
   1:    [ 3 . 2 . ]    (0,  3) (1,  4) (2)       [ 3 4 2 0 1 ]
   2:    [ 2 . 2 . ]    (0,  2) (1,  4) (3)       [ 2 4 0 3 1 ]
   3:    [ 1 . 2 . ]    (0,  1) (2,  4) (3)       [ 1 0 4 3 2 ]
   4:    [ . . 2 . ]    (0) (1) (2,  4) (3)       [ 0 1 4 3 2 ]
   5:    [ . . 1 . ]    (0) (1) (2,  3) (4)       [ 0 1 3 2 4 ]
   6:    [ 1 . 1 . ]    (0,  1) (2,  3) (4)       [ 1 0 3 2 4 ]
   7:    [ 2 . 1 . ]    (0,  2) (1,  3) (4)       [ 2 3 0 1 4 ]
   8:    [ 3 . 1 . ]    (0,  3) (1,  2) (4)       [ 3 2 1 0 4 ]
   9:    [ 4 . 1 . ]    (0,  4) (1,  2) (3)       [ 4 2 1 3 0 ]
  10:    [ 4 . . . ]    (0,  4) (1) (2) (3)       [ 4 1 2 3 0 ]
  11:    [ 3 . . . ]    (0,  3) (1) (2) (4)       [ 3 1 2 0 4 ]
  12:    [ 2 . . . ]    (0,  2) (1) (3) (4)       [ 2 1 0 3 4 ]
  13:    [ 1 . . . ]    (0,  1) (2) (3) (4)       [ 1 0 2 3 4 ]
  14:    [ . . . . ]    (0) (1) (2) (3) (4)       [ 0 1 2 3 4 ]
  15:    [ . 1 . . ]    (0) (1,  2) (3) (4)       [ 0 2 1 3 4 ]
  16:    [ . 2 . . ]    (0) (1,  3) (2) (4)       [ 0 3 2 1 4 ]
  17:    [ . 3 . . ]    (0) (1,  4) (2) (3)       [ 0 4 2 3 1 ]
  18:    [ . 3 . 1 ]    (0) (1,  4) (2,  3)       [ 0 4 3 2 1 ]
  19:    [ . 2 . 1 ]    (0) (1,  3) (2,  4)       [ 0 3 4 1 2 ]
  20:    [ . 1 . 1 ]    (0) (1,  2) (3,  4)       [ 0 2 1 4 3 ]
  21:    [ . . . 1 ]    (0) (1) (2) (3,  4)       [ 0 1 2 4 3 ]
  22:    [ 1 . . 1 ]    (0,  1) (2) (3,  4)       [ 1 0 2 4 3 ]
  23:    [ 2 . . 1 ]    (0,  2) (1) (3,  4)       [ 2 1 0 4 3 ]
  24:    [ 3 . . 1 ]    (0,  3) (1) (2,  4)       [ 3 1 4 0 2 ]
  25:    [ 4 . . 1 ]    (0,  4) (1) (2,  3)       [ 4 1 3 2 0 ]
\end{verbatim}
}
\caption{\label{fig:naf-invol}
Gray code for the length-4 non-adjacent forms in falling factorial base (left),
together with the corresponding involutions in cycle form (middle) and
in array form (right).
}
\end{figure}
%

The method of using a list of nonzero positions in the words to obtain loopless
methods for both \texttt{next()} and \texttt{prev()} has been implemented for
words where no two adjacent digits are nonzero (non-adjacent forms, NAF).

In the following algorithm for the computation of the successor
let $a_i$ be the digits of the length-$n$ NAF where $0\leq{}i<n$.
\begin{alg}[Next-NAF-SL]
Compute the successor of a non-adjacent form in subset-lex order.
\end{alg}
\begin{enumerate}
\item Let $j$ be the position of the last nonzero digit.
\item If $a_j$ is not maximal, increment it and return.
\item If $j+2 < n$, set $a_{j+2}=1$ (append new digit) and return.
\item Set $a_j=0$ (set last nonzero digit to zero).
\item If $j+1<n$, set $a_{j+1}=1$ (move digit right) and return.
\item If there is just one nonzero digit, stop.
\item Let $k$ be the position of the nearest nonzero digit left of $j$.
\item Set $a_k=a_k-1$ (decrement digit to the left).
\item If $a_k=0$, set $a_{k+1}=1$ (move digit right) and return.
\item Otherwise, $a_{k+2}=1$ (append digit two positions to the right).
\end{enumerate}
In the following implementation the array \texttt{iset[]}
holds the positions of the nonzero digits.
We only read the last or second last element,
the adjustments of the length (variable \texttt{ni})
have been left out in the description above.
The implementation handles all $n\geq{}0$ correctly,
for the all-zero word \texttt{iset[]} is of length one
and its only element points to the leftmost digit.
{\codesize
\begin{listing}{1}
// FILE:  src/comb/mixedradix-naf-subset-lex.h
class mixedradix_naf_subset_lex
{
    ulong *iset_;  // Set of positions of nonzero digits
    ulong *a_;    // digits
    ulong *m1_;   // nines (radix minus one) for each digit
    ulong ni_;    // number of elements in iset[]
    ulong n_;     // number of digits

    void first()
    {
        for (ulong k=0; k<n_; ++k)  a_[k] = 0;

        // iset[] initially with one element zero:
        iset_[0] = 0;
        ni_ = 1;

        if ( n_==0 )  // make things work for n == 0
        {
            m1_[0] = 0;
            a_[0] = 0;
        }
    }
\end{listing}
}%
The computation of the successor is
{\codesize
\begin{listing}{1}
    bool next()
    {
        ulong j = iset_[ni_-1];
        const ulong aj = a_[j] + 1;
        if ( aj <= m1_[j] )  // can increment last digit
        {
            a_[j] = aj;
            return true;
        }

        if ( j + 2 < n_ )  // can append new digit
        {
            iset_[ni_] = j + 2;
            a_[j+2] = 1;  // assume all m1[] are nonzero
            ++ni_;
            return true;
        }

        a_[j] = 0;  // set last nonzero digit to zero

        if ( j + 1 < n_ )  // can move last digit to the right
        {
            a_[j+1] = 1;
            iset_[ni_-1] = j + 1;
            return true;
        }

        if ( ni_ == 1 )  return false;  // current is last

        // Now we look to the left:
        const ulong k = iset_[ni_-2];  // nearest nonzero digit to the left
        const ulong ak = a_[k] - 1;    // decrement digit to the left
        a_[k] = ak;
        if ( ak == 0 )  // move digit one to the right
        {
            a_[k+1] = 1;
            iset_[ni_-2] = k + 1;
            --ni_;
            return true;
        }
        else  // append digit two positions to the right
        {
            a_[k+2] = 1;
            iset_[ni_-1] = k + 2;
            return true;
        }
    }
\end{listing}
}%
The update via \texttt{next()} takes about 7 cycles.
The updates for co-lexicographic order and the Gray code
respectively take about 15 and 21 cycles.

We note that the falling factorial NAFs of length $n$ are one-to-one
with involutions (self-inverse permutations) of $n+1$ elements.
Process the NAF from left to right; for $a_i=0$ let the
next unused element of the permutation be a fixed point
(and mark it as used); for $a_i\neq{}0$ put the next unused
element in a cycle with the $a_i$th unused element (and mark both as used).
A (simplistic) implementation of this method
is given in the program
\fileref{demo/comb/perm-involution-naf-demo.cc}.

%
\begin{figure}
{\jjbls
\begin{verbatim}
        colex        Gray code      subset-lex
   0:  ........      .1..1..1      1.......  =  { 0 }
   1:  .......1      .1..1...      1.1.....  =  { 0, 2 }
   2:  ......1.      .1..1.1.      1.1.1...  =  { 0, 2, 4 }
   3:  .....1..      .1....1.      1.1.1.1.  =  { 0, 2, 4, 6 }
   4:  .....1.1      .1......      1.1.1..1  =  { 0, 2, 4, 7 }
   5:  ....1...      .1.....1      1.1..1..  =  { 0, 2, 5 }
   6:  ....1..1      .1...1.1      1.1..1.1  =  { 0, 2, 5, 7 }
   7:  ....1.1.      .1...1..      1.1...1.  =  { 0, 2, 6 }
   8:  ...1....      .1.1.1..      1.1....1  =  { 0, 2, 7 }
   9:  ...1...1      .1.1.1.1      1..1....  =  { 0, 3 }
  10:  ...1..1.      .1.1...1      1..1.1..  =  { 0, 3, 5 }
  11:  ...1.1..      .1.1....      1..1.1.1  =  { 0, 3, 5, 7 }
  12:  ...1.1.1      .1.1..1.      1..1..1.  =  { 0, 3, 6 }
  13:  ..1.....      ...1..1.      1..1...1  =  { 0, 3, 7 }
  14:  ..1....1      ...1....      1...1...  =  { 0, 4 }
  15:  ..1...1.      ...1...1      1...1.1.  =  { 0, 4, 6 }
  16:  ..1..1..      ...1.1.1      1...1..1  =  { 0, 4, 7 }
  17:  ..1..1.1      ...1.1..      1....1..  =  { 0, 5 }
  18:  ..1.1...      .....1..      1....1.1  =  { 0, 5, 7 }
  19:  ..1.1..1      .....1.1      1.....1.  =  { 0, 6 }
  20:  ..1.1.1.      .......1      1......1  =  { 0, 7 }
  21:  .1......      ........      .1......  =  { 1 }
  22:  .1.....1      ......1.      .1.1....  =  { 1, 3 }
  23:  .1....1.      ....1.1.      .1.1.1..  =  { 1, 3, 5 }
  24:  .1...1..      ....1...      .1.1.1.1  =  { 1, 3, 5, 7 }
  25:  .1...1.1      ....1..1      .1.1..1.  =  { 1, 3, 6 }
  26:  .1..1...      ..1.1..1      .1.1...1  =  { 1, 3, 7 }
  27:  .1..1..1      ..1.1...      .1..1...  =  { 1, 4 }
  28:  .1..1.1.      ..1.1.1.      .1..1.1.  =  { 1, 4, 6 }
  29:  .1.1....      ..1...1.      .1..1..1  =  { 1, 4, 7 }
  30:  .1.1...1      ..1.....      .1...1..  =  { 1, 5 }
  31:  .1.1..1.      ..1....1      .1...1.1  =  { 1, 5, 7 }
  32:  .1.1.1..      ..1..1.1      .1....1.  =  { 1, 6 }
  33:  .1.1.1.1      ..1..1..      .1.....1  =  { 1, 7 }
  34:  1.......      1.1..1..      ..1.....  =  { 2 }
  35:  1......1      1.1..1.1      ..1.1...  =  { 2, 4 }
  36:  1.....1.      1.1....1      ..1.1.1.  =  { 2, 4, 6 }
  37:  1....1..      1.1.....      ..1.1..1  =  { 2, 4, 7 }
  38:  1....1.1      1.1...1.      ..1..1..  =  { 2, 5 }
  39:  1...1...      1.1.1.1.      ..1..1.1  =  { 2, 5, 7 }
  40:  1...1..1      1.1.1...      ..1...1.  =  { 2, 6 }
  41:  1...1.1.      1.1.1..1      ..1....1  =  { 2, 7 }
  42:  1..1....      1...1..1      ...1....  =  { 3 }
  43:  1..1...1      1...1...      ...1.1..  =  { 3, 5 }
  44:  1..1..1.      1...1.1.      ...1.1.1  =  { 3, 5, 7 }
  45:  1..1.1..      1.....1.      ...1..1.  =  { 3, 6 }
  46:  1..1.1.1      1.......      ...1...1  =  { 3, 7 }
  47:  1.1.....      1......1      ....1...  =  { 4 }
  48:  1.1....1      1....1.1      ....1.1.  =  { 4, 6 }
  49:  1.1...1.      1....1..      ....1..1  =  { 4, 7 }
  50:  1.1..1..      1..1.1..      .....1..  =  { 5 }
  51:  1.1..1.1      1..1.1.1      .....1.1  =  { 5, 7 }
  52:  1.1.1...      1..1...1      ......1.  =  { 6 }
  53:  1.1.1..1      1..1....      .......1  =  { 7 }
  54:  1.1.1.1.      1..1..1.
\end{verbatim}
}
\caption{\label{fig:bitfibrep}
Binary non-adjacent forms of length 8
in co-lexicographic order, minimal-change order (Gray code),
and subset-lex order (with the corresponding sets without
consecutive elements).
}
\end{figure}
%

The binary NAFs of length 8 are shown in figure \ref{fig:bitfibrep}.
Algorithms for the generation of these NAFs as binary words
for both lexicographic order and the Gray code are given
in \cite[p.75-77]{fxtbook}.
Here we give the implementations for computing successor
and predecessor for the nonempty NAFs in subset-lex order.


The routine for the successor needs to
start with the word that has a single set bit at
the highest position of the desired word length.
{\codesize
\begin{listing}{1}
// FILE:  src/bits/fibrep-subset-lexrev.h
ulong next_subset_lexrev_fib(ulong x)
{
    ulong x0 = x & -x;  // lowest bit
    ulong xs = x0 >> 2;
    if ( xs != 0 )  // easy case: set bit right of lowest bit
    {
        x |= xs;
        return  x;
    }
    else  // lowest bit at index 0 or 1
    {
        if ( x0 == 2 )  // at index 1
        {
            x -= 1;
            return x;
        }

        x ^= x0;  // clear lowest bit
        x0 = x & -x;  // new lowest bit ...
        x0 >>= 1;  x -= x0;  // ... is moved one to the right
        return  x;
    }
}
\end{listing}
}%
The all-zero word is returned as successor
of the word whose value is 1.
The routine for the predecessor can
be started with the all-zero word.
{\codesize
\begin{listing}{1}
ulong prev_subset_lexrev_fib(ulong x)
{
    ulong x0 = x & -x;  // lowest bit
    if ( x & (x0<<2) )  // easy case: next higher bit is set
    {
        x ^= x0;  // clear lowest bit
        return x;
    }
    else
    {
        x += x0;       // move lowest bit to the left and
        x |= (x0!=1);  // set rightmost bit unless blocked by next bit
        return x;
    }
}
\end{listing}
}%
%

\section{Compositions}

%
\begin{figure}
{\jjbls
\begin{verbatim}
       0:  ......    [ 1 1 1 1 1 1 1 ]
       1:  .....1    [ 1 1 1 1 1 2 ]
       2:  ....1.    [ 1 1 1 1 2 1 ]
       3:  ....11    [ 1 1 1 1 3 ]
       4:  ...1..    [ 1 1 1 2 1 1 ]
       5:  ...1.1    [ 1 1 1 2 2 ]
       6:  ...11.    [ 1 1 1 3 1 ]
       7:  ...111    [ 1 1 1 4 ]
       8:  ..1...    [ 1 1 2 1 1 1 ]
       9:  ..1..1    [ 1 1 2 1 2 ]
      10:  ..1.1.    [ 1 1 2 2 1 ]
      11:  ..1.11    [ 1 1 2 3 ]
      12:  ..11..    [ 1 1 3 1 1 ]
      13:  ..11.1    [ 1 1 3 2 ]
      14:  ..111.    [ 1 1 4 1 ]
      15:  ..1111    [ 1 1 5 ]
      16:  .1....    [ 1 2 1 1 1 1 ]
      17:  .1...1    [ 1 2 1 1 2 ]
      18:  .1..1.    [ 1 2 1 2 1 ]
      19:  .1..11    [ 1 2 1 3 ]
      20:  .1.1..    [ 1 2 2 1 1 ]
      21:  .1.1.1    [ 1 2 2 2 ]
      22:  .1.11.    [ 1 2 3 1 ]
      23:  .1.111    [ 1 2 4 ]
      24:  .11...    [ 1 3 1 1 1 ]
      25:  .11..1    [ 1 3 1 2 ]
      26:  .11.1.    [ 1 3 2 1 ]
      27:  .11.11    [ 1 3 3 ]
      28:  .111..    [ 1 4 1 1 ]
      29:  .111.1    [ 1 4 2 ]
      30:  .1111.    [ 1 5 1 ]
      31:  .11111    [ 1 6 ]
      32:  1.....    [ 2 1 1 1 1 1 ]
      33:  1....1    [ 2 1 1 1 2 ]
      34:  1...1.    [ 2 1 1 2 1 ]
      35:  1...11    [ 2 1 1 3 ]
      36:  1..1..    [ 2 1 2 1 1 ]
      37:  1..1.1    [ 2 1 2 2 ]
      38:  1..11.    [ 2 1 3 1 ]
      39:  1..111    [ 2 1 4 ]
      40:  1.1...    [ 2 2 1 1 1 ]
      41:  1.1..1    [ 2 2 1 2 ]
      42:  1.1.1.    [ 2 2 2 1 ]
      43:  1.1.11    [ 2 2 3 ]
      44:  1.11..    [ 2 3 1 1 ]
      45:  1.11.1    [ 2 3 2 ]
      46:  1.111.    [ 2 4 1 ]
      47:  1.1111    [ 2 5 ]
      48:  11....    [ 3 1 1 1 1 ]
      49:  11...1    [ 3 1 1 2 ]
      50:  11..1.    [ 3 1 2 1 ]
      51:  11..11    [ 3 1 3 ]
      52:  11.1..    [ 3 2 1 1 ]
      53:  11.1.1    [ 3 2 2 ]
      54:  11.11.    [ 3 3 1 ]
      55:  11.111    [ 3 4 ]
      56:  111...    [ 4 1 1 1 ]
      57:  111..1    [ 4 1 2 ]
      58:  111.1.    [ 4 2 1 ]
      59:  111.11    [ 4 3 ]
      60:  1111..    [ 5 1 1 ]
      61:  1111.1    [ 5 2 ]
      62:  11111.    [ 6 1 ]
      63:  111111    [ 7 ]
\end{verbatim}
}%
\caption{\label{fig:comp-nz}
The compositions of 7 together with their run-length
encodings as binary words (where dots denote zeros),
lexicographic order.
A succession of $k$ ones, followed by a zero, in the run-length encoding
stands for a part $k+1$, one trailing zero is implied.
}
\end{figure}
%

One of the most simple algorithms in combinatorial generation
may be the computation of the successor of a composition
represented as a list of parts for the lexicographic ordering.
In the following let $z$ be the last element of the composition.
\begin{samepage}
\begin{alg}[Next-Comp]
Compute the successor of a composition in lexicographic order.
\end{alg}
\begin{enumerate}
\item If there is just one part, stop (this is the last composition).
\item Add $1$ to the second last part (and remove $z$) and
 append $z-1$ ones at the end.
\end{enumerate}
\end{samepage}
\begin{samepage}
\begin{alg}[Prev-Comp]
Compute the predecessor of a composition in lexicographic order.
\end{alg}
\begin{enumerate}
\item If the number of parts is maximal (composition into all ones), stop.
\item If $z>1$, replace $z$ by $z-1,\,1$ (move one unit right); return
\item Otherwise, replace the tail $y,\, 1^q$ ($y$ followed by $q$ ones) by $y-1, q+1$.
\end{enumerate}
\end{samepage}
Figure \ref{fig:comp-nz} shows the compositions of 7
in lexicographic order.  The corresponding run-length encodings
appear in lexicographic order as well.

The algorithm for the successor can be made loopless
when care is taken that only ones are left beyond
the end of the current composition.

{\codesize
\begin{listing}{1}
// FILE:  src/comb/composition-nz.h
class composition_nz
// Compositions of n into positive parts, lexicographic order.
{
    ulong *a_;  // composition: a[1] + a[2] + ... + a[m] = n
    ulong n_;   // composition of n
    ulong m_;   // current composition is into m parts

    ulong next()
    // Return number of parts of generated composition.
    // Return zero if the current is the last composition.
    {
        if ( m_ <= 1 )  return 0;  // current is last

        // [*, Y, Z] --> [*, Y+1, 1, 1, 1, ..., 1]  (Z-1 trailing ones)
        a_[m_-1] += 1;
        const ulong z = a_[m_];
        a_[m_] = 1;
        // all parts a[m+1], a[m+2], ..., a[n] are already ==1
        m_ += z - 2;

        return  m_;
    }
\end{listing}
}%
%

%
\begin{figure}
{\jjbls
\begin{verbatim}
   0:  111111  [ 7 ]                      ......  [ 7 ]
   1:  11111.  [ 6 1 ]                    1.....  [ 1 6 ]
   2:  1111..  [ 5 1 1 ]                  11....  [ 1 1 5 ]
   3:  1111.1  [ 5 2 ]                    111...  [ 1 1 1 4 ]
   4:  111..1  [ 4 1 2 ]                  1111..  [ 1 1 1 1 3 ]
   5:  111...  [ 4 1 1 1 ]                11111.  [ 1 1 1 1 1 2 ]
   6:  111.1.  [ 4 2 1 ]                  111111  [ 1 1 1 1 1 1 1 ]
   7:  111.11  [ 4 3 ]                    1111.1  [ 1 1 1 1 2 1 ]
   8:  11..11  [ 3 1 3 ]                  111.1.  [ 1 1 1 2 2 ]
   9:  11..1.  [ 3 1 2 1 ]                111.11  [ 1 1 1 2 1 1 ]
  10:  11....  [ 3 1 1 1 1 ]              111..1  [ 1 1 1 3 1 ]
  11:  11...1  [ 3 1 1 2 ]                11.1..  [ 1 1 2 3 ]
  12:  11.1.1  [ 3 2 2 ]                  11.11.  [ 1 1 2 1 2 ]
  13:  11.1..  [ 3 2 1 1 ]                11.111  [ 1 1 2 1 1 1 ]
  14:  11.11.  [ 3 3 1 ]                  11.1.1  [ 1 1 2 2 1 ]
  15:  11.111  [ 3 4 ]                    11..1.  [ 1 1 3 2 ]
  16:  1..111  [ 2 1 4 ]                  11..11  [ 1 1 3 1 1 ]
  17:  1..11.  [ 2 1 3 1 ]                11...1  [ 1 1 4 1 ]
  18:  1..1..  [ 2 1 2 1 1 ]              1.1...  [ 1 2 4 ]
  19:  1..1.1  [ 2 1 2 2 ]                1.11..  [ 1 2 1 3 ]
  20:  1....1  [ 2 1 1 1 2 ]              1.111.  [ 1 2 1 1 2 ]
  21:  1.....  [ 2 1 1 1 1 1 ]            1.1111  [ 1 2 1 1 1 1 ]
  22:  1...1.  [ 2 1 1 2 1 ]              1.11.1  [ 1 2 1 2 1 ]
  23:  1...11  [ 2 1 1 3 ]                1.1.1.  [ 1 2 2 2 ]
  24:  1.1.11  [ 2 2 3 ]                  1.1.11  [ 1 2 2 1 1 ]
  25:  1.1.1.  [ 2 2 2 1 ]                1.1..1  [ 1 2 3 1 ]
  26:  1.1...  [ 2 2 1 1 1 ]              1..1..  [ 1 3 3 ]
  27:  1.1..1  [ 2 2 1 2 ]                1..11.  [ 1 3 1 2 ]
  28:  1.11.1  [ 2 3 2 ]                  1..111  [ 1 3 1 1 1 ]
  29:  1.11..  [ 2 3 1 1 ]                1..1.1  [ 1 3 2 1 ]
  30:  1.111.  [ 2 4 1 ]                  1...1.  [ 1 4 2 ]
  31:  1.1111  [ 2 5 ]                    1...11  [ 1 4 1 1 ]
  32:  ..1111  [ 1 1 5 ]                  1....1  [ 1 5 1 ]
  33:  ..111.  [ 1 1 4 1 ]                .1....  [ 2 5 ]
  34:  ..11..  [ 1 1 3 1 1 ]              .11...  [ 2 1 4 ]
  35:  ..11.1  [ 1 1 3 2 ]                .111..  [ 2 1 1 3 ]
  36:  ..1..1  [ 1 1 2 1 2 ]              .1111.  [ 2 1 1 1 2 ]
  37:  ..1...  [ 1 1 2 1 1 1 ]            .11111  [ 2 1 1 1 1 1 ]
  38:  ..1.1.  [ 1 1 2 2 1 ]              .111.1  [ 2 1 1 2 1 ]
  39:  ..1.11  [ 1 1 2 3 ]                .11.1.  [ 2 1 2 2 ]
  40:  ....11  [ 1 1 1 1 3 ]              .11.11  [ 2 1 2 1 1 ]
  41:  ....1.  [ 1 1 1 1 2 1 ]            .11..1  [ 2 1 3 1 ]
  42:  ......  [ 1 1 1 1 1 1 1 ]          .1.1..  [ 2 2 3 ]
  43:  .....1  [ 1 1 1 1 1 2 ]            .1.11.  [ 2 2 1 2 ]
  44:  ...1.1  [ 1 1 1 2 2 ]              .1.111  [ 2 2 1 1 1 ]
  45:  ...1..  [ 1 1 1 2 1 1 ]            .1.1.1  [ 2 2 2 1 ]
  46:  ...11.  [ 1 1 1 3 1 ]              .1..1.  [ 2 3 2 ]
  47:  ...111  [ 1 1 1 4 ]                .1..11  [ 2 3 1 1 ]
  48:  .1.111  [ 1 2 4 ]                  .1...1  [ 2 4 1 ]
  49:  .1.11.  [ 1 2 3 1 ]                ..1...  [ 3 4 ]
  50:  .1.1..  [ 1 2 2 1 1 ]              ..11..  [ 3 1 3 ]
  51:  .1.1.1  [ 1 2 2 2 ]                ..111.  [ 3 1 1 2 ]
  52:  .1...1  [ 1 2 1 1 2 ]              ..1111  [ 3 1 1 1 1 ]
  53:  .1....  [ 1 2 1 1 1 1 ]            ..11.1  [ 3 1 2 1 ]
  54:  .1..1.  [ 1 2 1 2 1 ]              ..1.1.  [ 3 2 2 ]
  55:  .1..11  [ 1 2 1 3 ]                ..1.11  [ 3 2 1 1 ]
  56:  .11.11  [ 1 3 3 ]                  ..1..1  [ 3 3 1 ]
  57:  .11.1.  [ 1 3 2 1 ]                ...1..  [ 4 3 ]
  58:  .11...  [ 1 3 1 1 1 ]              ...11.  [ 4 1 2 ]
  59:  .11..1  [ 1 3 1 2 ]                ...111  [ 4 1 1 1 ]
  60:  .111.1  [ 1 4 2 ]                  ...1.1  [ 4 2 1 ]
  61:  .111..  [ 1 4 1 1 ]                ....1.  [ 5 2 ]
  62:  .1111.  [ 1 5 1 ]                  ....11  [ 5 1 1 ]
  63:  .11111  [ 1 6 ]                    .....1  [ 6 1 ]
\end{verbatim}
}%
\caption{\label{fig:comp-nz-rl-and-sl}
The compositions of 7 together with their run-length
encodings as binary words (where dots denote zeros),
in an order corresponding to the complemented binary Gray code
(left) and in subset-lex order (right).
A succession of $k$ ones, followed by a zero, in the run-length encoding
stands for a part $k+1$, one trailing zero is implied (left).
The roles of ones and zeros are reversed for the subset-lex order (right).
}
\end{figure}
%

Figure \ref{fig:comp-nz-rl-and-sl} shows two orderings for the compositions.
The ordering corresponding to the (complemented) binary Gray code is shown in
the left columns.  We will call this order \jjterm{RL-order}, as the
compositions correspond to the run-lengths of the binary words in lexicographic
order.

For comparison with the new algorithm we give the (loopless) algorithms.
In the following let $m$ be the number of parts in the composition
and $x,y,z$ the last three parts.
\begin{samepage}
\begin{alg}[Next-Comp-RL]
Compute the successor of a composition in RL-order.
\end{alg}
\begin{enumerate}
\item If $m$ is odd:
 if $z\geq{}2$, replace $z$ by $z-1, 1$ and return;
 otherwise ($z=1$) replace $y,1$ by $y+1$ and return.
\item If $m$ is even:
 if $y\geq{}2$, replace $y,z$ by $y-1, 1, z$ and return;
 otherwise ($y=1$) replace $x,1,z$ by $x+1,z$ and return.
\end{enumerate}
\end{samepage}
The next algorithm is obtained from the previous simply
by swapping \eqq{even} and \eqq{odd} in the description.
\begin{samepage}
\begin{alg}[Prev-Comp-RL]
Compute the predecessor of a composition in RL-order.
\end{alg}
\begin{enumerate}
\item If $m$ is even:
 if $z\geq{}2$, replace $z$ by $z-1, 1$ and return;
 otherwise ($z=1$) replace $y,1$ by $y+1$ and return.
\item If $m$ is odd:
 if $y\geq{}2$, replace $y,z$ by $y-1, 1, z$ and return;
 otherwise ($y=1$) replace $x,1,z$ by $x+1,z$ and return.
\end{enumerate}
\end{samepage}
With each transition, at most three parts (at the end of
the composition) are changed.
The number of parts changes by $1$
with each step\footnote{See
 \fileref{src/comb/composition-nz-rl.h} for an implementation.}.

An alternative algorithm for this ordering is given
in \cite[ex.12, sect.7.2.1.1, p.308]{DEK4}, see
also \cite{misra}.

Now we give the algorithms for subset-lex order
(the first is essentially given in \cite{mansour-ll-comp}).

\begin{alg}[Next-Comp-SL]
Compute the successor of a composition in subset-lex order.
\end{alg}
\begin{enumerate}
\item If $z=1$ and there are at most two parts, stop.
\item If $z\geq{}2$, replace $z$ by $1, z-1$ (move all but one unit to the right); return.
\item Otherwise ($z=1$) add $1$ to the third last part and
  remove $z$ (move one unit two places left: $x,y,1$ is replaced by $x+1, y$).
\end{enumerate}
For the next algorithm let $y,\,z$ be the two last elements.
\begin{samepage}
\begin{alg}[Prev-Comp-SL]
Compute the predecessor of a composition in subset-lex order.
\end{alg}
\begin{enumerate}
\item If there is just one part, stop.
\item If $y=1$, replace $y,\,z$ by $z+1$ (add $z$ to the left); return.
\item Otherwise ($y\geq{}2$) replace $y,\,z$ by $y-1,\,z,\,1$ (move one unit two places right).
\end{enumerate}
\end{samepage}

At most two parts are changed by either method and
these span at most those three parts at the end of the composition.
Again, the number of parts changes by $1$ with each step.

Note that also the initializations are loopless
for the preceding two orderings.

We give the crucial parts of the implementation.
{\codesize
\begin{listing}{1}
// FILE: src/comb/composition-nz-subset-lex.h
class composition_nz_subset_lex
{
    ulong *a_;  // composition: a[1] + a[2] + ... + a[m] = n
    ulong n_;   // composition of n
    ulong m_;   // current composition is into m parts

    void first()
    {
        a_[0] = 0;
        a_[1] = n_;
        m_ = ( n_ ? 1 : 0 );
    }

    void last()
    {
        if ( n_ >= 2 )
        {
            a_[1] = n_ - 1;
            a_[2] = 1;
            for (ulong j=2; j<=n_; ++j)  a_[j] = 1;
            m_ = 2;
        }
        else
        {
            a_[1] = n_;
            m_ = n_;
        }
    }
\end{listing}
}%

For the methods \texttt{next()} and \texttt{prev()}
we give one of the two implementations found in the file.
Both methods return the number of parts in the generated
composition and return zero if there are no more
compositions.
{\codesize
\begin{listing}{1}
    ulong next()
    {
        const ulong z = a_[m_];
        if ( z<=1 )  // move one unit two places left
        { //   [*, X, Y, 1] --> [*, X+1, Y]
            if ( m_ <= 2 )  return 0;  // current is last
            m_ -= 1;
            a_[m_-1] += 1;
            return  m_;
        }
        else  // move all but one unit right
        { //   [*, Y, Z] --> [*, Y, 1, Z-1]
            a_[m_] = 1;
            m_ += 1;
            a_[m_] = z - 1;
            return  m_;
        }
    }
\end{listing}
}%
{\codesize
\begin{listing}{1}
    ulong prev()
    {
        if ( m_ <= 1 )  return 0;  // current is first

        const ulong y = a_[m_-1];
        if ( y==1 )  // add Z to left place
        { //   [*, 1, Z] --> [*, Z+1]
            const ulong z = a_[m_];
            a_[m_-1] = z + 1;
            a_[m_] = 1;
            m_ -= 1;
            return  m_;
        }
        else  // move one unit two places right
        { //   [*, Y, Z] --> [*, Y-1, Z, 1]
            a_[m_-1] = y - 1;
            m_ += 1;
            // a[m] == 1 already
            return  m_;
        }
    }
\end{listing}
}%

The method \texttt{next()} takes about
5 cycles for lexicographic order,
9 cycles for RL-order,
and 6 cycles for subset-lex order.
The method \texttt{prev()} takes about
10.5 cycles for lexicographic order
and is identical in performance
to \texttt{next()} for the other orders.

Note that the last parts in the successive compositions in
lexicographic order give the (one-based)
ruler function (sequence \jjseqref{A001511} in \cite{oeis}), see
\fileref{src/comb/ruler-func1.h} for the trivial implementation, compare to
\fileref{src/comb/ruler-func.h} (giving sequence \jjseqref{A007814} in \cite{oeis})
which uses the techniques from \cite{ehrlichll} and \cite{bitnergray}
as described in \cite[Algorithm L, p.290]{DEK4}.

An loopless implementation for the generation of the compositions into odd parts in
subset-lex order is given in \fileref{src/comb/composition-nz-odd-subset-lex.h}.
%

\subsubsection*{Ranking and unranking}
An unranking algorithm for the subset-lex order
is obtained by observing (see figure \ref{fig:comp-nz-rl-and-sl})
that the composition into one part has rank 0
and otherwise the first part and the remaining parts
are easily determined.
\begin{alg}[Unrank-Comp-SL]\label{alg:unrank-comp-sl}
Recursive routine $U( r,\, n,\, C[\,],\, m )$ for the computation
of the composition $C[\,]$ of $n$ with rank $r$ in subset-lex order.
The auxiliary variable $m\geq{}0$ is the index
of the part about to be written in $C[\,]$.
The initial call is $F(r,\, n,\, C[\,],\, 0)$.
\end{alg}
\begin{enumerate}
\item If $r=0$, set $C[m]=n$ and return $m+1$ (the number of parts).
\item Set $f=0$ (the tentative first part).
\item\label{alg:loop}
 Set $f=f+1$ and $t=2^{n-f-1}$.
\item If $r<t$ (can use part $f$),
 set $C[m]=f$ and return $F( r,\, n-f,\, C[\,],\, m+1 )$.
\item Set $r=r-t$ and go to step \ref{alg:loop}.
\end{enumerate}
We give an iterative implementation.
%
{\codesize
\begin{listing}{1}
// FILE:  src/comb/composition-nz-rank.cc
ulong composition_nz_subset_lex_unrank(ulong r, ulong *x, ulong n)
{
    if ( n==0 )  return 0;

    ulong m = 0;
    while ( true )
    {
        if ( r==0 )  // composition into one part
        {
            x[m++] = n;
            return m;
        }
        r -= 1;
        ulong t = 1UL << (n-1);
        for (ulong f=1; f<n; ++f)  // find first part f >= 1
        {
            t >>= 1;  // == 2**(n-f-1)
            // == number of compositions of n with first part f

            if ( r < t )  // first part is f
            {
                x[m++] = f;
                n -= f;
                break;
            }
            r -= t;
        }
    }
}
\end{listing}
}%

The following algorithm for computing the rank
is modeled as inverse of the method above.
\begin{alg}[Rank-Comp-SL]\label{alg:rank-comp-sl}
Computation of the rank $r$ of a composition $C[\,]$ of $n$ in subset-lex order.
\end{alg}
\begin{enumerate}
\item Set $r=0$ and $e=0$ (position of part under consideration).
\item\label{alg:loop2}
 Set $f=C[e]$ (part under consideration).
\item If $f=n$, return $r$.
\item Set $r=r+1$, $t=2^{n-1}$, and $n=n-f$.
\item While $f>1$, set $t=t/2$, $r=r+t$, and $f=f-1$.
\item Set $e=e+1$ and go to step \ref{alg:loop2}
\end{enumerate}
An implementation is
{\codesize
\begin{listing}{1}
ulong composition_nz_subset_lex_rank(const ulong *x, ulong m, ulong n)
{
    ulong r = 0;  // rank
    ulong e = 0;  // position of first part
    while ( e < m )
    {
        ulong f = x[e];
        if ( f==n )  return r;
        r += 1;
        ulong t = 1UL << (n-1);
        n -= f;
        while ( f > 1 )
        {
            t >>= 1;
            r += t;
            f -= 1;
        }
        e += 1;
    }
    return r;  // return r==0 for the empty composition
}
\end{listing}
}%
%

%
%
%

Conversion functions between binary words in lexicographic order and subset-lex
order are shown in \cite[pp.71-72]{fxtbook}.
{\codesize
\begin{listing}{1}
// FILE: src/bits/bitlex.h
ulong negidx2lexrev(ulong k)
{
    ulong z = 0;
    ulong h = highest_one(k);
    while ( k )
    {
        while ( 0==(h&k) )  h >>= 1;
        z ^= h;
        ++k;
        k &= h - 1;
    }

    return  z;
}

ulong lexrev2negidx(ulong x)
{
    if ( 0==x )  return 0;
    ulong h = x & -x;  // lowest one
    ulong r = (h-1);
    while ( x^=h )
    {
        r += (h-1);
        h = x & -x;  // next higher one
    }
    r += h;  // highest bit
    return  r;
}
\end{listing}
}%
Based on these, alternative ranking and unranking functions
can be given\footnote{See
 \fileref{src/comb/composition-nz-rank.cc}
where the corresponding routines for all three orders shown here are implemented.}.

Ranking and unranking methods for the compositions
into odd parts can be obtained by replacing
$2^{n-f-1}$ (number of compositions of $n$ with first part $f$)
in the algorithms \ref{alg:unrank-comp-sl} and \ref{alg:rank-comp-sl}
by the Fibonacci numbers $F_{n-f-1}$
(number of compositions of $n$ into odd parts with first part $f$),
where $n\geq{}1$, $f<n$, and $f$ odd
(see sequence \jjseqref{A242086} in \cite{oeis}).
%

\section{Partitions as weakly increasing lists of parts}

We now give algorithms for the computation
of the successor for partitions represented
as  weakly increasing lists of parts.
The generation of all partitions in this representation
is the subject of \cite{kelleher-thesis}.

\subsection{All partitions}

%
\begin{figure}
{\jjbls
\begin{verbatim}
   1:   [ 1 1 1 1 1 1 1 1 1 1 1 ]      [ 11 ]
   2:   [ 1 1 1 1 1 1 1 1 1 2 ]        [ 1 10 ]
   3:   [ 1 1 1 1 1 1 1 1 3 ]          [ 1 1 9 ]
   4:   [ 1 1 1 1 1 1 1 2 2 ]          [ 1 1 1 8 ]
   5:   [ 1 1 1 1 1 1 1 4 ]            [ 1 1 1 1 7 ]
   6:   [ 1 1 1 1 1 1 2 3 ]            [ 1 1 1 1 1 6 ]
   7:   [ 1 1 1 1 1 1 5 ]              [ 1 1 1 1 1 1 5 ]
   8:   [ 1 1 1 1 1 2 2 2 ]            [ 1 1 1 1 1 1 1 4 ]
   9:   [ 1 1 1 1 1 2 4 ]              [ 1 1 1 1 1 1 1 1 3 ]
  10:   [ 1 1 1 1 1 3 3 ]              [ 1 1 1 1 1 1 1 1 1 2 ]
  11:   [ 1 1 1 1 1 6 ]                [ 1 1 1 1 1 1 1 1 1 1 1 ]
  12:   [ 1 1 1 1 2 2 3 ]              [ 1 1 1 1 1 1 1 2 2 ]
  13:   [ 1 1 1 1 2 5 ]                [ 1 1 1 1 1 1 2 3 ]
  14:   [ 1 1 1 1 3 4 ]                [ 1 1 1 1 1 2 4 ]
  15:   [ 1 1 1 1 7 ]                  [ 1 1 1 1 1 2 2 2 ]
  16:   [ 1 1 1 2 2 2 2 ]              [ 1 1 1 1 1 3 3 ]
  17:   [ 1 1 1 2 2 4 ]                [ 1 1 1 1 2 5 ]
  18:   [ 1 1 1 2 3 3 ]                [ 1 1 1 1 2 2 3 ]
  19:   [ 1 1 1 2 6 ]                  [ 1 1 1 1 3 4 ]
  20:   [ 1 1 1 3 5 ]                  [ 1 1 1 2 6 ]
  21:   [ 1 1 1 4 4 ]                  [ 1 1 1 2 2 4 ]
  22:   [ 1 1 1 8 ]                    [ 1 1 1 2 2 2 2 ]
  23:   [ 1 1 2 2 2 3 ]                [ 1 1 1 2 3 3 ]
  24:   [ 1 1 2 2 5 ]                  [ 1 1 1 3 5 ]
  25:   [ 1 1 2 3 4 ]                  [ 1 1 1 4 4 ]
  26:   [ 1 1 2 7 ]                    [ 1 1 2 7 ]
  27:   [ 1 1 3 3 3 ]                  [ 1 1 2 2 5 ]
  28:   [ 1 1 3 6 ]                    [ 1 1 2 2 2 3 ]
  29:   [ 1 1 4 5 ]                    [ 1 1 2 3 4 ]
  30:   [ 1 1 9 ]                      [ 1 1 3 6 ]
  31:   [ 1 2 2 2 2 2 ]                [ 1 1 3 3 3 ]
  32:   [ 1 2 2 2 4 ]                  [ 1 1 4 5 ]
  33:   [ 1 2 2 3 3 ]                  [ 1 2 8 ]
  34:   [ 1 2 2 6 ]                    [ 1 2 2 6 ]
  35:   [ 1 2 3 5 ]                    [ 1 2 2 2 4 ]
  36:   [ 1 2 4 4 ]                    [ 1 2 2 2 2 2 ]
  37:   [ 1 2 8 ]                      [ 1 2 2 3 3 ]
  38:   [ 1 3 3 4 ]                    [ 1 2 3 5 ]
  39:   [ 1 3 7 ]                      [ 1 2 4 4 ]
  40:   [ 1 4 6 ]                      [ 1 3 7 ]
  41:   [ 1 5 5 ]                      [ 1 3 3 4 ]
  42:   [ 1 10 ]                       [ 1 4 6 ]
  43:   [ 2 2 2 2 3 ]                  [ 1 5 5 ]
  44:   [ 2 2 2 5 ]                    [ 2 9 ]
  45:   [ 2 2 3 4 ]                    [ 2 2 7 ]
  46:   [ 2 2 7 ]                      [ 2 2 2 5 ]
  47:   [ 2 3 3 3 ]                    [ 2 2 2 2 3 ]
  48:   [ 2 3 6 ]                      [ 2 2 3 4 ]
  49:   [ 2 4 5 ]                      [ 2 3 6 ]
  50:   [ 2 9 ]                        [ 2 3 3 3 ]
  51:   [ 3 3 5 ]                      [ 2 4 5 ]
  52:   [ 3 4 4 ]                      [ 3 8 ]
  53:   [ 3 8 ]                        [ 3 3 5 ]
  54:   [ 4 7 ]                        [ 3 4 4 ]
  55:   [ 5 6 ]                        [ 4 7 ]
  56:   [ 11 ]                         [ 5 6 ]
\end{verbatim}
}
\caption{\label{fig:partition-asc}
The partitions of 11 as weakly increasing lists of parts,
lexicographic order (left) and subset-lex order (right).
}
\end{figure}
%

Figure \ref{fig:partition-asc}
shows the partitions of 11 as weakly increasing lists of parts,
in lexicographic order (left) and in subset-lex order (right).
We first give a description of the computation of
the successor in lexicographic order.
The three last parts of the partition are denoted by $x$, $y$, and $z$.
\begin{samepage}
\begin{alg}[Next-Part-Asc]
Compute the successor of a partition in lexicographic order.
\end{alg}
\begin{enumerate}
\item If there is just one part, stop.
\item If $z-1 < y+1$, replace $y,z$ by $y+z$; return.
\item Otherwise, change $y$ to $y+1$ and append
 parts $y+1$ as long as there are at least $y+1$ units left.
\item Add the remaining units to the last part.
\end{enumerate}
\end{samepage}
The second step of the algorithm can be merged
into the other steps, as it is a special case of them.
{\codesize
\begin{listing}{1}
// FILE:  src/comb/partition-asc.h
class partition_asc
{
    ulong *a_;  // partition: a[1] + a[2] + ... + a[m] = n
    ulong n_;   // integer partitions of n
    ulong m_;   // current partition has m parts
\end{listing}
}%
The implementation is correct for all $n\geq{}0$,
for $n=0$ the empty list of parts is generated.
{\codesize
\begin{listing}{1}
    ulong next()
    // Return number of parts of generated partition.
    // Return zero if the current partition is the last.
    {
        if ( m_ <= 1 )  return 0;  // current is last

        ulong z1 = a_[m_] - 1;  // take one unit from last part
        m_ -= 1;
        const ulong y1 = a_[m_] + 1;  // add one unit to previous part

        while ( y1 <= z1 )  // can put part Y+1
        {
            a_[m_] = y1;
            z1 -= y1;
            m_ += 1;
        }
        a_[m_] = y1 + z1;  // add remaining units to last part

        return  m_;
    }
\end{listing}
}%
%

The computation of the successor in subset-lex order is loopless.
\begin{samepage}
\begin{alg}[Next-Part-Asc-SL]
Compute the successor of a partition in subset-lex order.
A sentinel zero shall precede list of elements.
\end{alg}
\begin{enumerate}
\item If $z\geq{}2\,y$, replace $y,z$ by $y,y,z-y$ (extend to the right); return.
\item If $z-1\geq{}y+1$, replace $y,z$ by $y+1,z-1$ (add one unit to the left); return.
\item If $z=1$ (the all-ones partition) do the following.
  Stop if the number of parts is $\leq{}3$,
  otherwise replace the tail $[1,1,1,1]$ by $[2,2]$ and return.
\item If the number of parts is 2, stop.
\item Replace $x, y, z$ by $x+1, y+z-1$ (add one unit to second left, add rest to end) and return.
\end{enumerate}
\end{samepage}
Note that with each update all but the last two parts
are the same as in the predecessor.
This can be an advantage over the lexicographic order for
computations where partial results for prefixes can be reused.

{\codesize
\begin{listing}{1}
// FILE: src/comb/partition-asc-subset-lex.h
class partition_asc_subset_lex
{
    ulong *a_;  // partition: a[1] + a[2] + ... + a[m] = n
    ulong n_;   // partition of n
    ulong m_;   // current partition has m parts

    explicit partition_asc_subset_lex(ulong n)
    {
        n_ = n;
        a_ = new ulong[n_+1+(n_==0)];
        // sentinel a[0] set in first()
        first();
    }

    void first()
    {
        a_[0] = 1;  // sentinel: read (once) by test for right-extension
        a_[1] = n_ + (n_==0);  // use partitions of n=1 for n=0 internally
        m_ = 1;
    }
\end{listing}
}%
The implementation is again correct for all $n\geq{}0$.
{\codesize
\begin{listing}{1}
    ulong next()
    // Loopless algorithm.
    {
        ulong y = a_[m_-1];  // may read sentinel a[0]
        ulong z = a_[m_];

        if ( z >= 2*y )  // extend to the right:
        {  //   [*, Y, Z] --> [*, Y, Y, Z-Y]
            a_[m_] = y;
            a_[m_+1] = z - y;  // >= y
            ++m_;
            return  m_;
        }

        z -= 1;  y += 1;
        if ( z >= y )  // add one unit to the left:
        {  //   [*, Y, Z] --> [*, Y+1, Z-1]
            a_[m_-1] = y;
            a_[m_] = z;
            return  m_;
        }

        if ( z == 0 )  // all-ones partition
        {
            if ( n_ <= 3 )  return 0;  // current is last

            //   [1, ..., 1, 1, 1, 1] -->  [1, ..., 2, 2]
            m_ -= 2;
            a_[m_] = 2;  a_[m_-1] = 2;
            return m_;
        }

        // add one unit to second left, add rest to end:
        //   [*, X, Y, Z] --> [*, X+1, Y+Z-1]
        a_[m_-2] += 1;
        a_[m_-1] += z;
        m_ -= 1;
        m_ -= (m_ == 1);  // last if partition is into one part
        return  m_;
    }
\end{listing}
}%

The updates via \texttt{next()}
take about 15 cycles for lexicographic order
and about 13 cycles for subset-lex order.

\subsection{Partitions into odd and into distinct parts}

%
\begin{figure}
{\jjbls
\begin{verbatim}
   0:  [ 19 ]              [ 19 ]
   1:  [ 1 18 ]            [ 1 1 17 ]
   2:  [ 1 2 16 ]          [ 1 1 1 1 15 ]
   3:  [ 1 2 3 13 ]        [ 1 1 1 1 1 1 13 ]
   4:  [ 1 2 3 4 9 ]       [ 1 1 1 1 1 1 1 1 11 ]
   5:  [ 1 2 3 5 8 ]       [ 1 1 1 1 1 1 1 1 1 1 9 ]
   6:  [ 1 2 3 6 7 ]       [ 1 1 1 1 1 1 1 1 1 1 1 1 7 ]
   7:  [ 1 2 4 12 ]        [ 1 1 1 1 1 1 1 1 1 1 1 1 1 1 5 ]
   8:  [ 1 2 4 5 7 ]       [ 1 1 1 1 1 1 1 1 1 1 1 1 1 1 1 1 3 ]
   9:  [ 1 2 5 11 ]        [ 1 1 1 1 1 1 1 1 1 1 1 1 1 1 1 1 1 1 1 ]
  10:  [ 1 2 6 10 ]        [ 1 1 1 1 1 1 1 1 1 1 1 1 1 3 3 ]
  11:  [ 1 2 7 9 ]         [ 1 1 1 1 1 1 1 1 1 1 1 3 5 ]
  12:  [ 1 3 15 ]          [ 1 1 1 1 1 1 1 1 1 1 3 3 3 ]
  13:  [ 1 3 4 11 ]        [ 1 1 1 1 1 1 1 1 1 3 7 ]
  14:  [ 1 3 4 5 6 ]       [ 1 1 1 1 1 1 1 1 1 5 5 ]
  15:  [ 1 3 5 10 ]        [ 1 1 1 1 1 1 1 1 3 3 5 ]
  16:  [ 1 3 6 9 ]         [ 1 1 1 1 1 1 1 3 9 ]
  17:  [ 1 3 7 8 ]         [ 1 1 1 1 1 1 1 3 3 3 3 ]
  18:  [ 1 4 14 ]          [ 1 1 1 1 1 1 1 5 7 ]
  19:  [ 1 4 5 9 ]         [ 1 1 1 1 1 1 3 3 7 ]
  20:  [ 1 4 6 8 ]         [ 1 1 1 1 1 1 3 5 5 ]
  21:  [ 1 5 13 ]          [ 1 1 1 1 1 3 11 ]
  22:  [ 1 5 6 7 ]         [ 1 1 1 1 1 3 3 3 5 ]
  23:  [ 1 6 12 ]          [ 1 1 1 1 1 5 9 ]
  24:  [ 1 7 11 ]          [ 1 1 1 1 1 7 7 ]
  25:  [ 1 8 10 ]          [ 1 1 1 1 3 3 9 ]
  26:  [ 2 17 ]            [ 1 1 1 1 3 3 3 3 3 ]
  27:  [ 2 3 14 ]          [ 1 1 1 1 3 5 7 ]
  28:  [ 2 3 4 10 ]        [ 1 1 1 1 5 5 5 ]
  29:  [ 2 3 5 9 ]         [ 1 1 1 3 13 ]
  30:  [ 2 3 6 8 ]         [ 1 1 1 3 3 3 7 ]
  31:  [ 2 4 13 ]          [ 1 1 1 3 3 5 5 ]
  32:  [ 2 4 5 8 ]         [ 1 1 1 5 11 ]
  33:  [ 2 4 6 7 ]         [ 1 1 1 7 9 ]
  34:  [ 2 5 12 ]          [ 1 1 3 3 11 ]
  35:  [ 2 6 11 ]          [ 1 1 3 3 3 3 5 ]
  36:  [ 2 7 10 ]          [ 1 1 3 5 9 ]
  37:  [ 2 8 9 ]           [ 1 1 3 7 7 ]
  38:  [ 3 16 ]            [ 1 1 5 5 7 ]
  39:  [ 3 4 12 ]          [ 1 3 15 ]
  40:  [ 3 4 5 7 ]         [ 1 3 3 3 9 ]
  41:  [ 3 5 11 ]          [ 1 3 3 3 3 3 3 ]
  42:  [ 3 6 10 ]          [ 1 3 3 5 7 ]
  43:  [ 3 7 9 ]           [ 1 3 5 5 5 ]
  44:  [ 4 15 ]            [ 1 5 13 ]
  45:  [ 4 5 10 ]          [ 1 7 11 ]
  46:  [ 4 6 9 ]           [ 1 9 9 ]
  47:  [ 4 7 8 ]           [ 3 3 13 ]
  48:  [ 5 14 ]            [ 3 3 3 3 7 ]
  49:  [ 5 6 8 ]           [ 3 3 3 5 5 ]
  50:  [ 6 13 ]            [ 3 5 11 ]
  51:  [ 7 12 ]            [ 3 7 9 ]
  52:  [ 8 11 ]            [ 5 5 9 ]
  53:  [ 9 10 ]            [ 5 7 7 ]
\end{verbatim}
}
\caption{\label{fig:partition-asc-odd-dist}
The partitions of 19 into distinct parts (left) and odd parts (right) in subset-lex order.}
\end{figure}
%

It is not very difficult to adapt the algorithms
to specialized partitions.
We give the algorithms for partitions into distinct and odd parts
by their implementations.

Computation of the successor for partitions into distinct parts
 in lexicographic order:
{\codesize
\begin{listing}{1}
// FILE:  src/comb/partition-dist-asc.h
    ulong next()
    {
        if ( m_ <= 1 )  return 0;  // current is last

        ulong s = a_[m_] + a_[m_-1];
        ulong k = a_[m_-1] + 1;
        m_ -= 1;
        // split s into k, k+1, k+2, ..., y, z  where z >= y + 1:
        while ( s >= ( k + (k+1) ) )
        {
            a_[m_] = k;
            s -= k;
            k += 1;
            m_ += 1;
        }

        a_[m_] = s;

        return  m_;
    }
\end{listing}
}%

Computation of the successor for partitions into odd parts
 in lexicographic order:
{\codesize
\begin{listing}{1}
// FILE:  src/comb/partition-odd-asc.h
    ulong next()
    {
        const ulong z = a_[m_];   // can read sentinel a[0] if n==0
        const ulong y = a_[m_-1]; // can read sentinel a[0] (a[-1] for n==0)
        ulong s = y + z;  // sum of parts we scan over

        ulong k;  // min value of next term
        if ( z >= y+4 )   // add last 2 terms
        {
            if ( m_ == 1 )  return 0;  // current is last
            k = y + 2;
            a_[m_-1] = k;
            s -= k;
        }
        else  // add last 3 terms
        {
            if ( m_ <= 2 )  return 0;  // current is last
            const ulong x = a_[m_-2];
            s += x;
            k = x + 2;
            m_ -= 2;
        }

        const ulong k2 = k + k;
        while ( s >= k2 + k )
        {
            a_[m_] = k;  s -= k;  m_ += 1;
            a_[m_] = k;  s -= k;  m_ += 1;
        }

        a_[m_] = s;
        return  m_;
    }
\end{listing}
}%

All routines in this section return the number of parts in the generated
partition and return zero if there are no more partitions.


Computation of the successor for partitions into distinct parts
 in subset-lex order:
{\codesize
\begin{listing}{1}
// FILE:  src/comb/partition-dist-asc-subset-lex.h
    ulong next()
    // Loopless algorithm.
    {
        ulong y = a_[m_-1];  // may read sentinel a[0]
        ulong z = a_[m_];

        if ( z >= 2*y + 3 )  // can extend to the right
        {  // [*, Y, Z] --> [*, Y, Y+1, Z-1]
            y += 1;
            a_[m_] = y;
            a_[m_+1] = z - y;  // >= y
            ++m_;
            return  m_;
        }
        else  // add to the left
        {
            z -= 1;  y += 1;

            if ( z > y )  // add one unit to the left
            {  // [*, Y, Z] --> [*, Y+1, Z-1]

                if ( m_<=1 )  return 0;  // current is last

                a_[m_-1] = y;
                a_[m_] = z;
                return  m_;
            }
            else  // add to one unit second left
                  // and combine last with second last
            {  // [*, X, Y, Z] --> [*, X+1, Y+Z]

                if ( m_<=2 )  return 0;  // current is last

                a_[m_-2] += 1;
                a_[m_-1] += z;
                --m_;
                return  m_;
            }
        }
    }
\end{listing}
}%
Computation of the successor for partitions into odd parts
 in subset-lex order:
{\codesize
\begin{listing}{1}
// FILE: src/comb/partitions-odd-asc-subset-lex.h
    ulong next()
    // Loopless algorithm.
    {
        ulong y = a_[m_-1];  // may read sentinel a[0]
        ulong z = a_[m_];

        if ( z >= 3*y )  // can extend to the right
        {  // [*, Y, Z] --> [*, Y, Y, Y, Z-2*Y]
            a_[m_] = y;
            a_[m_+1] = y;
            a_[m_+2] = z - 2 * y;
            m_ += 2;
            return  m_;
        }

        if ( m_ >= n_ )  // all-ones partition
        {
            if ( n_ <= 5 )  return 0;  // current is last

            //   [1, ..., 1, 1, 1, 1, 1, 1] -->  [1, ..., 3, 3]
            m_ -= 4;
            a_[m_] = 3;
            a_[m_-1] = 3;
            return m_;
        }

        ulong z2 = z - 2;
        ulong y2 = y + 2;

        if ( z2 >= y2 )  // add 2 units to the left
        {  // [*, Y, Z] --> [*, Y+2, Z-2]
            a_[m_-1] = y2;
            a_[m_] = z2;
            return  m_;
        }

        if ( m_==2 )  // current is last (happens only for n even)
            return 0;

        // here m >= 3;  add to second or third left:

        ulong x2 = a_[m_-2] + 2;
        ulong s = z + y - 2;

        if ( x2 <= z2 )
        // add 2 units to third left, repeat part, put rest at end
        {  // [*, X, Y, Z] --> [*, X+2, X+2, Y+Z-2 -X-2]
            a_[m_-2] = x2;
            a_[m_-1] = x2;
            a_[m_] = s - x2;
            return  m_;
        }

        if ( m_==3 )  // current is last (happens only for n odd)
            return 0;

        // add 2 units to third left, combine rest into second left
        // [*, W, X, Y, Z] --> [*, W+2, X+Y+Z-2]
        a_[m_-3] += 2;
        a_[m_-2] += s;
        m_ -= 2;
        return  m_;
    }
\end{listing}
}%

The updates via \texttt{next()} of partitions into distinct parts take about 11
cycles for lexicographic order and about 9 cycles for subset-lex order.  The
respective figures for the partitions into odd parts are 13 and 13.5 cycles.

Algorithms for ranking and unranking are obtained by replacing $2^{n-f-1}$ (the
number of compositions of $n$ with first part $f$) in algorithms
\ref{alg:unrank-comp-sl} and \ref{alg:rank-comp-sl} by the expression for the
number of partitions of $n$ (of the desired kind) with first part $f$.

%
%
%
%
%
%
%
%
%
%
%

\clearpage
\section{Subset-lex order for restricted growth strings (RGS)}

By modifying the algorithms for the successor and predecessor (\ref{alg:next-sl}
and \ref{alg:prev-sl}) for subset-lex order to adhere to certain conditions,
one obtains the equivalents for RGS.
These can often be given without any difficulty.
We give two examples,
RGS for set partitions and
RGS for $k$-ary Dyck words.

\subsection{RGS for set partitions}

%
\begin{figure}
{\jjbls
\begin{verbatim}
      lexicographic                            subset-lex
 1:  [ . . . . . ]  {1 2 3 4 5}           [ . . . . . ]  {1 2 3 4 5}
 2:  [ . . . . 1 ]  {1 2 3 4} {5}         [ . 1 . . . ]  {1 3 4 5} {2}
 3:  [ . . . 1 . ]  {1 2 3 5} {4}         [ . 1 1 . . ]  {1 4 5} {2 3}
 4:  [ . . . 1 1 ]  {1 2 3} {4 5}         [ . 1 2 . . ]  {1 4 5} {2} {3}
 5:  [ . . . 1 2 ]  {1 2 3} {4} {5}       [ . 1 2 1 . ]  {1 5} {2 4} {3}
 6:  [ . . 1 . . ]  {1 2 4 5} {3}         [ . 1 2 2 . ]  {1 5} {2} {3 4}
 7:  [ . . 1 . 1 ]  {1 2 4} {3 5}         [ . 1 2 3 . ]  {1 5} {2} {3} {4}
 8:  [ . . 1 . 2 ]  {1 2 4} {3} {5}       [ . 1 2 3 1 ]  {1} {2 5} {3} {4}
 9:  [ . . 1 1 . ]  {1 2 5} {3 4}         [ . 1 2 3 2 ]  {1} {2} {3 5} {4}
10:  [ . . 1 1 1 ]  {1 2} {3 4 5}         [ . 1 2 3 3 ]  {1} {2} {3} {4 5}
11:  [ . . 1 1 2 ]  {1 2} {3 4} {5}       [ . 1 2 3 4 ]  {1} {2} {3} {4} {5}
12:  [ . . 1 2 . ]  {1 2 5} {3} {4}       [ . 1 2 2 1 ]  {1} {2 5} {3 4}
13:  [ . . 1 2 1 ]  {1 2} {3 5} {4}       [ . 1 2 2 2 ]  {1} {2} {3 4 5}
14:  [ . . 1 2 2 ]  {1 2} {3} {4 5}       [ . 1 2 2 3 ]  {1} {2} {3 4} {5}
15:  [ . . 1 2 3 ]  {1 2} {3} {4} {5}     [ . 1 2 1 1 ]  {1} {2 4 5} {3}
16:  [ . 1 . . . ]  {1 3 4 5} {2}         [ . 1 2 1 2 ]  {1} {2 4} {3 5}
17:  [ . 1 . . 1 ]  {1 3 4} {2 5}         [ . 1 2 1 3 ]  {1} {2 4} {3} {5}
18:  [ . 1 . . 2 ]  {1 3 4} {2} {5}       [ . 1 2 . 1 ]  {1 4} {2 5} {3}
19:  [ . 1 . 1 . ]  {1 3 5} {2 4}         [ . 1 2 . 2 ]  {1 4} {2} {3 5}
20:  [ . 1 . 1 1 ]  {1 3} {2 4 5}         [ . 1 2 . 3 ]  {1 4} {2} {3} {5}
21:  [ . 1 . 1 2 ]  {1 3} {2 4} {5}       [ . 1 1 1 . ]  {1 5} {2 3 4}
22:  [ . 1 . 2 . ]  {1 3 5} {2} {4}       [ . 1 1 2 . ]  {1 5} {2 3} {4}
23:  [ . 1 . 2 1 ]  {1 3} {2 5} {4}       [ . 1 1 2 1 ]  {1} {2 3 5} {4}
24:  [ . 1 . 2 2 ]  {1 3} {2} {4 5}       [ . 1 1 2 2 ]  {1} {2 3} {4 5}
25:  [ . 1 . 2 3 ]  {1 3} {2} {4} {5}     [ . 1 1 2 3 ]  {1} {2 3} {4} {5}
26:  [ . 1 1 . . ]  {1 4 5} {2 3}         [ . 1 1 1 1 ]  {1} {2 3 4 5}
27:  [ . 1 1 . 1 ]  {1 4} {2 3 5}         [ . 1 1 1 2 ]  {1} {2 3 4} {5}
28:  [ . 1 1 . 2 ]  {1 4} {2 3} {5}       [ . 1 1 . 1 ]  {1 4} {2 3 5}
29:  [ . 1 1 1 . ]  {1 5} {2 3 4}         [ . 1 1 . 2 ]  {1 4} {2 3} {5}
30:  [ . 1 1 1 1 ]  {1} {2 3 4 5}         [ . 1 . 1 . ]  {1 3 5} {2 4}
31:  [ . 1 1 1 2 ]  {1} {2 3 4} {5}       [ . 1 . 2 . ]  {1 3 5} {2} {4}
32:  [ . 1 1 2 . ]  {1 5} {2 3} {4}       [ . 1 . 2 1 ]  {1 3} {2 5} {4}
33:  [ . 1 1 2 1 ]  {1} {2 3 5} {4}       [ . 1 . 2 2 ]  {1 3} {2} {4 5}
34:  [ . 1 1 2 2 ]  {1} {2 3} {4 5}       [ . 1 . 2 3 ]  {1 3} {2} {4} {5}
35:  [ . 1 1 2 3 ]  {1} {2 3} {4} {5}     [ . 1 . 1 1 ]  {1 3} {2 4 5}
36:  [ . 1 2 . . ]  {1 4 5} {2} {3}       [ . 1 . 1 2 ]  {1 3} {2 4} {5}
37:  [ . 1 2 . 1 ]  {1 4} {2 5} {3}       [ . 1 . . 1 ]  {1 3 4} {2 5}
38:  [ . 1 2 . 2 ]  {1 4} {2} {3 5}       [ . 1 . . 2 ]  {1 3 4} {2} {5}
39:  [ . 1 2 . 3 ]  {1 4} {2} {3} {5}     [ . . 1 . . ]  {1 2 4 5} {3}
40:  [ . 1 2 1 . ]  {1 5} {2 4} {3}       [ . . 1 1 . ]  {1 2 5} {3 4}
41:  [ . 1 2 1 1 ]  {1} {2 4 5} {3}       [ . . 1 2 . ]  {1 2 5} {3} {4}
42:  [ . 1 2 1 2 ]  {1} {2 4} {3 5}       [ . . 1 2 1 ]  {1 2} {3 5} {4}
43:  [ . 1 2 1 3 ]  {1} {2 4} {3} {5}     [ . . 1 2 2 ]  {1 2} {3} {4 5}
44:  [ . 1 2 2 . ]  {1 5} {2} {3 4}       [ . . 1 2 3 ]  {1 2} {3} {4} {5}
45:  [ . 1 2 2 1 ]  {1} {2 5} {3 4}       [ . . 1 1 1 ]  {1 2} {3 4 5}
46:  [ . 1 2 2 2 ]  {1} {2} {3 4 5}       [ . . 1 1 2 ]  {1 2} {3 4} {5}
47:  [ . 1 2 2 3 ]  {1} {2} {3 4} {5}     [ . . 1 . 1 ]  {1 2 4} {3 5}
48:  [ . 1 2 3 . ]  {1 5} {2} {3} {4}     [ . . 1 . 2 ]  {1 2 4} {3} {5}
49:  [ . 1 2 3 1 ]  {1} {2 5} {3} {4}     [ . . . 1 . ]  {1 2 3 5} {4}
50:  [ . 1 2 3 2 ]  {1} {2} {3 5} {4}     [ . . . 1 1 ]  {1 2 3} {4 5}
51:  [ . 1 2 3 3 ]  {1} {2} {3} {4 5}     [ . . . 1 2 ]  {1 2 3} {4} {5}
52:  [ . 1 2 3 4 ]  {1} {2} {3} {4} {5}   [ . . . . 1 ]  {1 2 3 4} {5}
\end{verbatim}
}
\caption{\label{fig:setpart-rgs}
Restricted growth strings and corresponding set partitions in
lexicographic and subset-lex order.
}
\end{figure}
%

%
We consider the restricted growth strings (RGS)
$a_0,\,a_1,\,\ldots,\,a_{n-1}$ such that
$a_0=0$ and $a_k \leq{} 1 + \max(a_0,\, a_1,\, \ldots,\, a_{k-1})$.
These RGS are counted by the Bell numbers, see sequence \jjseqref{A000110} in \cite{oeis}.
Figure \ref{fig:setpart-rgs} shows the all such RGS of length 5
in lexicographic and subset-lex order.
The set partitions are obtained by putting
all $i$ such that $a_i=k$ into the same set.

In the following algorithm we assume that
 $m_k=1 + \max(a_0,\, a_1,\, \ldots,\, a_{k-1})$
and that there is a sentinel $a_{-1}=1$.
\begin{alg}[Next-Setpart-RGS]
Compute the successor in subset-lex order.
\end{alg}
\begin{enumerate}
\item Let $j$ be the index of the last nonzero digit ($j=1$ for the all-zero RGS).
\item If $a_j < m_j$, set $a_j=a_j+1$ and return.
\item If $j+1 < n$, set $a_{j+1}=1$ and return.
\item Set $a_j=0$.
\item Set $j$ to the position of the next nonzero digit to the left.
 If $j=-1$, stop.
\item Set $a_j=a_j-1$ and $a_{j+1}=1$.
\end{enumerate}
The implementation has to take care of updating
the array $m[]$ which has an additional (write-only)
element at its end.
{\codesize
\begin{listing}{1}
// FILE:  src/comb/setpart-rgs-subset-lex.h
class setpart_rgs_subset_lex
{
    ulong *a_;  // digits of the RGS
    ulong *m_;  // maximum + 1 in prefix
    ulong tr_;  // current track
    ulong n_;   // number of digits in RGS
\end{listing}
}%
The computation of the successor is correct for all $n\geq{}0$,
we omit the constructor, which sets $m_0=1$ for $n=0$ to
cover this case.
{\codesize
\begin{listing}{1}
    bool next()
    {
        ulong j = tr_;
        if ( a_[j] < m_[j] )   // easy case 1: can increment track
        {
            if ( n_ <= 1 )  return false;  // handle n <= 1 correctly
            a_[j] += 1;
            return true;
        }

        const ulong j1 = j + 1;
        if ( j1 < n_ )  // easy case 2: can attach
        {
            m_[j1] = m_[j] + 1;
            a_[j1] = +1;
            tr_ = j1;
            return true;
        }

        a_[j] = 0;
        m_[j] = m_[j-1];

        // Find nonzero track to the left:
        do  { --j; }  while ( a_[j] == 0 );  // can read sentinel

        if ( (long)j < 0 )  return false;  // current is last

        if ( a_[j] == m_[j] )  m_[j+1] = m_[j];
        a_[j] -= 1;

        ++j;
        a_[j] = 1;
        tr_ = j;
        return true;
    }
\end{listing}
}%

An update takes about 9.5 cycles for subset-lex order
and 12.5 cycles for lexicographic order.

%
%
%

\subsection{RGS for $k$-ary Dyck words}

%
\begin{figure}
{\jjbls
\begin{verbatim}
   1:  [ . . . . ]    1..1..1..1..          [ . . . . ]    1..1..1..1..
   2:  [ . . . 1 ]    1..1..1.1...          [ . 1 . . ]    1.1...1..1..
   3:  [ . . . 2 ]    1..1..11....          [ . 2 . . ]    11....1..1..
   4:  [ . . 1 . ]    1..1.1...1..          [ . 2 1 . ]    11...1...1..
   5:  [ . . 1 1 ]    1..1.1..1...          [ . 2 2 . ]    11..1....1..
   6:  [ . . 1 2 ]    1..1.1.1....          [ . 2 3 . ]    11.1.....1..
   7:  [ . . 1 3 ]    1..1.11.....          [ . 2 4 . ]    111......1..
   8:  [ . . 2 . ]    1..11....1..          [ . 2 4 1 ]    111.....1...
   9:  [ . . 2 1 ]    1..11...1...          [ . 2 4 2 ]    111....1....
  10:  [ . . 2 2 ]    1..11..1....          [ . 2 4 3 ]    111...1.....
  11:  [ . . 2 3 ]    1..11.1.....          [ . 2 4 4 ]    111..1......
  12:  [ . . 2 4 ]    1..111......          [ . 2 4 5 ]    111.1.......
  13:  [ . 1 . . ]    1.1...1..1..          [ . 2 4 6 ]    1111........
  14:  [ . 1 . 1 ]    1.1...1.1...          [ . 2 3 1 ]    11.1....1...
  15:  [ . 1 . 2 ]    1.1...11....          [ . 2 3 2 ]    11.1...1....
  16:  [ . 1 1 . ]    1.1..1...1..          [ . 2 3 3 ]    11.1..1.....
  17:  [ . 1 1 1 ]    1.1..1..1...          [ . 2 3 4 ]    11.1.1......
  18:  [ . 1 1 2 ]    1.1..1.1....          [ . 2 3 5 ]    11.11.......
  19:  [ . 1 1 3 ]    1.1..11.....          [ . 2 2 1 ]    11..1...1...
  20:  [ . 1 2 . ]    1.1.1....1..          [ . 2 2 2 ]    11..1..1....
  21:  [ . 1 2 1 ]    1.1.1...1...          [ . 2 2 3 ]    11..1.1.....
  22:  [ . 1 2 2 ]    1.1.1..1....          [ . 2 2 4 ]    11..11......
  23:  [ . 1 2 3 ]    1.1.1.1.....          [ . 2 1 1 ]    11...1..1...
  24:  [ . 1 2 4 ]    1.1.11......          [ . 2 1 2 ]    11...1.1....
  25:  [ . 1 3 . ]    1.11.....1..          [ . 2 1 3 ]    11...11.....
  26:  [ . 1 3 1 ]    1.11....1...          [ . 2 . 1 ]    11....1.1...
  27:  [ . 1 3 2 ]    1.11...1....          [ . 2 . 2 ]    11....11....
  28:  [ . 1 3 3 ]    1.11..1.....          [ . 1 1 . ]    1.1..1...1..
  29:  [ . 1 3 4 ]    1.11.1......          [ . 1 2 . ]    1.1.1....1..
  30:  [ . 1 3 5 ]    1.111.......          [ . 1 3 . ]    1.11.....1..
  31:  [ . 2 . . ]    11....1..1..          [ . 1 3 1 ]    1.11....1...
  32:  [ . 2 . 1 ]    11....1.1...          [ . 1 3 2 ]    1.11...1....
  33:  [ . 2 . 2 ]    11....11....          [ . 1 3 3 ]    1.11..1.....
  34:  [ . 2 1 . ]    11...1...1..          [ . 1 3 4 ]    1.11.1......
  35:  [ . 2 1 1 ]    11...1..1...          [ . 1 3 5 ]    1.111.......
  36:  [ . 2 1 2 ]    11...1.1....          [ . 1 2 1 ]    1.1.1...1...
  37:  [ . 2 1 3 ]    11...11.....          [ . 1 2 2 ]    1.1.1..1....
  38:  [ . 2 2 . ]    11..1....1..          [ . 1 2 3 ]    1.1.1.1.....
  39:  [ . 2 2 1 ]    11..1...1...          [ . 1 2 4 ]    1.1.11......
  40:  [ . 2 2 2 ]    11..1..1....          [ . 1 1 1 ]    1.1..1..1...
  41:  [ . 2 2 3 ]    11..1.1.....          [ . 1 1 2 ]    1.1..1.1....
  42:  [ . 2 2 4 ]    11..11......          [ . 1 1 3 ]    1.1..11.....
  43:  [ . 2 3 . ]    11.1.....1..          [ . 1 . 1 ]    1.1...1.1...
  44:  [ . 2 3 1 ]    11.1....1...          [ . 1 . 2 ]    1.1...11....
  45:  [ . 2 3 2 ]    11.1...1....          [ . . 1 . ]    1..1.1...1..
  46:  [ . 2 3 3 ]    11.1..1.....          [ . . 2 . ]    1..11....1..
  47:  [ . 2 3 4 ]    11.1.1......          [ . . 2 1 ]    1..11...1...
  48:  [ . 2 3 5 ]    11.11.......          [ . . 2 2 ]    1..11..1....
  49:  [ . 2 4 . ]    111......1..          [ . . 2 3 ]    1..11.1.....
  50:  [ . 2 4 1 ]    111.....1...          [ . . 2 4 ]    1..111......
  51:  [ . 2 4 2 ]    111....1....          [ . . 1 1 ]    1..1.1..1...
  52:  [ . 2 4 3 ]    111...1.....          [ . . 1 2 ]    1..1.1.1....
  53:  [ . 2 4 4 ]    111..1......          [ . . 1 3 ]    1..1.11.....
  54:  [ . 2 4 5 ]    111.1.......          [ . . . 1 ]    1..1..1.1...
  55:  [ . 2 4 6 ]    1111........          [ . . . 2 ]    1..1..11....
\end{verbatim}
}
\caption{\label{fig:dyck-words}
The 55 RGS for the 3-ary Dyck words of length 12, in
lexicographic order (left) and subset-lex order (right).
}
\end{figure}
%

Figure \ref{fig:dyck-words} shows
the 55 RGS for the 3-ary Dyck words of length 12, in
both lexicographic and subset-lex order.
The $j$th value in each RGS contains the distance of the position $j$th one in
the Dyck word from its maximal value $k\cdot{}j$
The sequences of numbers of $k$-ary Dyck words are entries
\jjseqref{A000108} ($k=2$, Catalan numbers),
\jjseqref{A001764} ($k=3$),
\jjseqref{A002293} ($k=4$), and
\jjseqref{A002294} ($k=5$) in \cite{oeis}.

It shall suffice to give the implementation for the (loopless) computation of
the predecessor in subset-lex order.
A sentinel \texttt{a[-1]=+1} is used.
{\codesize
\begin{listing}{1}
// FILE: src/comb/dyck-rgs-subset-lex.h
class dyck_rgs_subset_lex
{
    ulong *a_;  // digits of the RGS: a_[k] <= as[k-1] + 1
    ulong tr_;  // current track
    ulong n_;   // number of digits in RGS
    ulong i_;   // k-ary Dyck words: i = k - 1

    void last()
    {
        for (ulong k=0; k<n_; ++k)  a_[k] = 0;
        tr_ = n_ - 1;
        // make things work for n <= 1:
        if ( n_==0 )
        {
            tr_ = 0;
            a_[0] = 1;
        }
        if ( n_>=2 )  a_[tr_] = i_;
    }
\end{listing}
}%
All cases $n\geq{}0$ are handled correctly.
{\codesize
\begin{listing}{1}
    bool prev()
    // Loopless algorithm.
    {
        if ( n_<=1 )  return false;  // just one RGS

        ulong j = tr_;
        if ( a_[j] > 1 )   // can decrement track
        {
            a_[j] -= 1;
            return true;
        }

        const ulong aj = a_[j];  // zero or one

        a_[j] = 0;
        --j;

        if ( a_[j] == a_[j-1] + i_ )  // move track to the left
        {
            --tr_;
            return true;
        }

        if ( j==0 )  // current or next is last
        {
            if ( aj == 0 )  return false;
            return true;
        }

        a_[j] += 1;  // increment left digit
        tr_ = n_ - 1;  // move to right end
        a_[tr_] = a_[tr_-1] + i_;  // set to max value
        return true;
    }
\end{listing}
}%
One update takes about $10$ cycles for both
lexicographic and subset-lex order.

For the generation of other restricted growth strings in subset-lex order,
see the following files.
{\jjbls
\begin{verbatim}
// Ascent sequences, see OEIS sequence A022493:
src/comb/ascent-rgs-subset-lex.h
src/comb/ascent-rgs.h  // lexicographic order

// RGS for Catalan objects, see OEIS sequence A000108:
src/comb/catalan-rgs-subset-lex.h  (loopless prev())
src/comb/catalan-rgs.h  // lexicographic order

// Standard Young tableaux, represented as ballot sequences,
// see OEIS sequence A000085:
src/comb/young-tab-rgs-subset-lex.h
src/comb/young-tab-rgs.h  // lexicographic order
\end{verbatim}
}%

\section{A variant of the subset-lex order}

%
\begin{figure}
\vspace*{-10mm}
{\jjbls
\begin{verbatim}
    0:  [ . . . . ]    {  }
    1:  [ . . . 1 ]    { 3 }
    2:  [ . . . 2 ]    { 3, 3 }
    3:  [ . . . 3 ]    { 3, 3, 3 }
    4:  [ . . 1 3 ]    { 3, 3, 3, 2 }
    5:  [ . . 2 3 ]    { 3, 3, 3, 2, 2 }
    6:  [ . 1 2 3 ]    { 3, 3, 3, 2, 2, 1 }
    7:  [ . 2 2 3 ]    { 3, 3, 3, 2, 2, 1, 1 }
    8:  [ 1 2 2 3 ]    { 3, 3, 3, 2, 2, 1, 1, 0 }
    9:  [ 1 1 2 3 ]    { 3, 3, 3, 2, 2, 1, 0 }
   10:  [ 1 . 2 3 ]    { 3, 3, 3, 2, 2, 0 }
   11:  [ . 1 1 3 ]    { 3, 3, 3, 2, 1 }
   12:  [ . 2 1 3 ]    { 3, 3, 3, 2, 1, 1 }
   13:  [ 1 2 1 3 ]    { 3, 3, 3, 2, 1, 1, 0 }
   14:  [ 1 1 1 3 ]    { 3, 3, 3, 2, 1, 0 }
   15:  [ 1 . 1 3 ]    { 3, 3, 3, 2, 0 }
   16:  [ . 1 . 3 ]    { 3, 3, 3, 1 }
   17:  [ . 2 . 3 ]    { 3, 3, 3, 1, 1 }
   18:  [ 1 2 . 3 ]    { 3, 3, 3, 1, 1, 0 }
   19:  [ 1 1 . 3 ]    { 3, 3, 3, 1, 0 }
   20:  [ 1 . . 3 ]    { 3, 3, 3, 0 }
   21:  [ . . 1 2 ]    { 3, 3, 2 }
   22:  [ . . 2 2 ]    { 3, 3, 2, 2 }
   23:  [ . 1 2 2 ]    { 3, 3, 2, 2, 1 }
   24:  [ . 2 2 2 ]    { 3, 3, 2, 2, 1, 1 }
   25:  [ 1 2 2 2 ]    { 3, 3, 2, 2, 1, 1, 0 }
   26:  [ 1 1 2 2 ]    { 3, 3, 2, 2, 1, 0 }
   27:  [ 1 . 2 2 ]    { 3, 3, 2, 2, 0 }
   28:  [ . 1 1 2 ]    { 3, 3, 2, 1 }
   29:  [ . 2 1 2 ]    { 3, 3, 2, 1, 1 }
   30:  [ 1 2 1 2 ]    { 3, 3, 2, 1, 1, 0 }
   31:  [ 1 1 1 2 ]    { 3, 3, 2, 1, 0 }
   32:  [ 1 . 1 2 ]    { 3, 3, 2, 0 }
   33:  [ . 1 . 2 ]    { 3, 3, 1 }
   34:  [ . 2 . 2 ]    { 3, 3, 1, 1 }
   35:  [ 1 2 . 2 ]    { 3, 3, 1, 1, 0 }
   36:  [ 1 1 . 2 ]    { 3, 3, 1, 0 }
   37:  [ 1 . . 2 ]    { 3, 3, 0 }
   38:  [ . . 1 1 ]    { 3, 2 }
   39:  [ . . 2 1 ]    { 3, 2, 2 }
   40:  [ . 1 2 1 ]    { 3, 2, 2, 1 }
   41:  [ . 2 2 1 ]    { 3, 2, 2, 1, 1 }
   42:  [ 1 2 2 1 ]    { 3, 2, 2, 1, 1, 0 }
   43:  [ 1 1 2 1 ]    { 3, 2, 2, 1, 0 }
   44:  [ 1 . 2 1 ]    { 3, 2, 2, 0 }
   45:  [ . 1 1 1 ]    { 3, 2, 1 }
   46:  [ . 2 1 1 ]    { 3, 2, 1, 1 }
   47:  [ 1 2 1 1 ]    { 3, 2, 1, 1, 0 }
   48:  [ 1 1 1 1 ]    { 3, 2, 1, 0 }
   49:  [ 1 . 1 1 ]    { 3, 2, 0 }
   50:  [ . 1 . 1 ]    { 3, 1 }
   51:  [ . 2 . 1 ]    { 3, 1, 1 }
   52:  [ 1 2 . 1 ]    { 3, 1, 1, 0 }
   53:  [ 1 1 . 1 ]    { 3, 1, 0 }
   54:  [ 1 . . 1 ]    { 3, 0 }
   55:  [ . . 1 . ]    { 2 }
   56:  [ . . 2 . ]    { 2, 2 }
   57:  [ . 1 2 . ]    { 2, 2, 1 }
   58:  [ . 2 2 . ]    { 2, 2, 1, 1 }
   59:  [ 1 2 2 . ]    { 2, 2, 1, 1, 0 }
   60:  [ 1 1 2 . ]    { 2, 2, 1, 0 }
   61:  [ 1 . 2 . ]    { 2, 2, 0 }
   62:  [ . 1 1 . ]    { 2, 1 }
   63:  [ . 2 1 . ]    { 2, 1, 1 }
   64:  [ 1 2 1 . ]    { 2, 1, 1, 0 }
   65:  [ 1 1 1 . ]    { 2, 1, 0 }
   66:  [ 1 . 1 . ]    { 2, 0 }
   67:  [ . 1 . . ]    { 1 }
   68:  [ . 2 . . ]    { 1, 1 }
   69:  [ 1 2 . . ]    { 1, 1, 0 }
   70:  [ 1 1 . . ]    { 1, 0 }
   71:  [ 1 . . . ]    { 0 }
\end{verbatim}
}
\caption{\label{fig:subset-lexrev}
Subsets of the set $\{0^1, 1^2, 2^2, 3^3 \}$ in subset-lexrev order.
Dots denote zeros in the (generalized) characteristic words.
Note the sets are printed starting with the largest element.
}
\end{figure}
%

An order obtained by processing the positions in the characteristic
words as for subset-lex order but with reversed priorities is
shown in figure \ref{fig:subset-lexrev}.
We will call this ordering \jjterm{subset-lexrev}.
The algorithms for computing successor and predecessor
are easily obtained from those for the subset-lex order,
we just give the implementation for \texttt{prev()},
which is (again) loopless.
{\codesize
\begin{listing}{1}
// FILE: src/comb/mixedradix-subset-lexrev.h
    bool prev()
    {
        ulong j = tr_;
        if ( a_[j] > 1 )  // easy case: just decrement
        {
            a_[j] -= 1;
            return true;
        }

        a_[j] = 0;
        ++j;  // now looking at next track to the right

        if ( j >= n_ )  // was on rightmost track (last two steps)
        {
            bool q = ( a_[j] != 0 );
            a_[j] = 0;
            return q;
        }

        if ( a_[j] == m1_[j] )  // semi-easy case: move track to left
        {
            tr_ = j;  // move track one right
            return true;
        }
        else
        {
            a_[j] += 1;  // increment digit to the right
            j = 0;
            a_[j] = m1_[j];  // set leftmost digit = nine
            tr_ = j;  // move to leftmost track
            return true;
        }
    }
\end{listing}
}%
%

%
\begin{figure}
{\jjbls
\begin{verbatim}
          lex            colex         subset-lexrev     subset-lex
 1:  [ . . . . . ]    [ . . . . . ]    [ . . . . . ]    [ . . . . . ]
 2:  [ . . . . 1 ]    [ . . . . 1 ]    [ . . . . 1 ]    [ . 1 2 3 3 ]
 3:  [ . . . . 2 ]    [ . . . 1 1 ]    [ . . . . 2 ]    [ . 1 2 3 4 ]
 4:  [ . . . . 3 ]    [ . . 1 1 1 ]    [ . . . . 3 ]    [ . 1 2 2 2 ]
 5:  [ . . . . 4 ]    [ . 1 1 1 1 ]    [ . . . . 4 ]    [ . 1 2 2 3 ]
 6:  [ . . . 1 1 ]    [ . . . . 2 ]    [ . . . 1 4 ]    [ . 1 2 2 4 ]
 7:  [ . . . 1 2 ]    [ . . . 1 2 ]    [ . . . 2 4 ]    [ . 1 1 3 3 ]
 8:  [ . . . 1 3 ]    [ . . 1 1 2 ]    [ . . . 3 4 ]    [ . 1 1 3 4 ]
 9:  [ . . . 1 4 ]    [ . 1 1 1 2 ]    [ . . 1 3 4 ]    [ . 1 1 2 2 ]
10:  [ . . . 2 2 ]    [ . . . 2 2 ]    [ . . 2 3 4 ]    [ . 1 1 2 3 ]
11:  [ . . . 2 3 ]    [ . . 1 2 2 ]    [ . 1 2 3 4 ]    [ . 1 1 2 4 ]
12:  [ . . . 2 4 ]    [ . 1 1 2 2 ]    [ . 1 1 3 4 ]    [ . 1 1 1 1 ]
13:  [ . . . 3 3 ]    [ . . 2 2 2 ]    [ . . 1 2 4 ]    [ . 1 1 1 2 ]
14:  [ . . . 3 4 ]    [ . 1 2 2 2 ]    [ . . 2 2 4 ]    [ . 1 1 1 3 ]
15:  [ . . 1 1 1 ]    [ . . . . 3 ]    [ . 1 2 2 4 ]    [ . 1 1 1 4 ]
16:  [ . . 1 1 2 ]    [ . . . 1 3 ]    [ . 1 1 2 4 ]    [ . . 2 3 3 ]
17:  [ . . 1 1 3 ]    [ . . 1 1 3 ]    [ . . 1 1 4 ]    [ . . 2 3 4 ]
18:  [ . . 1 1 4 ]    [ . 1 1 1 3 ]    [ . 1 1 1 4 ]    [ . . 2 2 2 ]
19:  [ . . 1 2 2 ]    [ . . . 2 3 ]    [ . . . 1 3 ]    [ . . 2 2 3 ]
20:  [ . . 1 2 3 ]    [ . . 1 2 3 ]    [ . . . 2 3 ]    [ . . 2 2 4 ]
21:  [ . . 1 2 4 ]    [ . 1 1 2 3 ]    [ . . . 3 3 ]    [ . . 1 3 3 ]
22:  [ . . 1 3 3 ]    [ . . 2 2 3 ]    [ . . 1 3 3 ]    [ . . 1 3 4 ]
23:  [ . . 1 3 4 ]    [ . 1 2 2 3 ]    [ . . 2 3 3 ]    [ . . 1 2 2 ]
24:  [ . . 2 2 2 ]    [ . . . 3 3 ]    [ . 1 2 3 3 ]    [ . . 1 2 3 ]
25:  [ . . 2 2 3 ]    [ . . 1 3 3 ]    [ . 1 1 3 3 ]    [ . . 1 2 4 ]
26:  [ . . 2 2 4 ]    [ . 1 1 3 3 ]    [ . . 1 2 3 ]    [ . . 1 1 1 ]
27:  [ . . 2 3 3 ]    [ . . 2 3 3 ]    [ . . 2 2 3 ]    [ . . 1 1 2 ]
28:  [ . . 2 3 4 ]    [ . 1 2 3 3 ]    [ . 1 2 2 3 ]    [ . . 1 1 3 ]
29:  [ . 1 1 1 1 ]    [ . . . . 4 ]    [ . 1 1 2 3 ]    [ . . 1 1 4 ]
30:  [ . 1 1 1 2 ]    [ . . . 1 4 ]    [ . . 1 1 3 ]    [ . . . 3 3 ]
31:  [ . 1 1 1 3 ]    [ . . 1 1 4 ]    [ . 1 1 1 3 ]    [ . . . 3 4 ]
32:  [ . 1 1 1 4 ]    [ . 1 1 1 4 ]    [ . . . 1 2 ]    [ . . . 2 2 ]
33:  [ . 1 1 2 2 ]    [ . . . 2 4 ]    [ . . . 2 2 ]    [ . . . 2 3 ]
34:  [ . 1 1 2 3 ]    [ . . 1 2 4 ]    [ . . 1 2 2 ]    [ . . . 2 4 ]
35:  [ . 1 1 2 4 ]    [ . 1 1 2 4 ]    [ . . 2 2 2 ]    [ . . . 1 1 ]
36:  [ . 1 1 3 3 ]    [ . . 2 2 4 ]    [ . 1 2 2 2 ]    [ . . . 1 2 ]
37:  [ . 1 1 3 4 ]    [ . 1 2 2 4 ]    [ . 1 1 2 2 ]    [ . . . 1 3 ]
38:  [ . 1 2 2 2 ]    [ . . . 3 4 ]    [ . . 1 1 2 ]    [ . . . 1 4 ]
39:  [ . 1 2 2 3 ]    [ . . 1 3 4 ]    [ . 1 1 1 2 ]    [ . . . . 1 ]
40:  [ . 1 2 2 4 ]    [ . 1 1 3 4 ]    [ . . . 1 1 ]    [ . . . . 2 ]
41:  [ . 1 2 3 3 ]    [ . . 2 3 4 ]    [ . . 1 1 1 ]    [ . . . . 3 ]
42:  [ . 1 2 3 4 ]    [ . 1 2 3 4 ]    [ . 1 1 1 1 ]    [ . . . . 4 ]
\end{verbatim}
}
\caption{\label{fig:catalan-step-rgs}
RGS corresponding to certain lattice paths (see text)
in lexicographic, co-lexicographic, subset-lexrev, and subset-lex order.
}
\end{figure}

One use of the subset-lexrev order is the development of algorithms for the
generation of orderings for objects where the subset-lex order exhibits no
simple structure.

As an example we consider the restricted growth strings corresponding
to lattice paths from $(0,0)$ to $(n,n)$ with steps
$(+1,0)$ and $(+1,+1)$ that do not go below the diagonal.
These RGS are words $a_0,\,a_1, \ldots,\, a_{n-1}$
such that $a_0=0$, $a_j\leq{}j$, and $a_{j+1}\geq{}a_j$
(the final $a_n=n$ is omitted).
These RGS are counted by the Catalan numbers,
see sequence \jjseqref{A000108} in \cite{oeis}.

Figure \ref{fig:catalan-step-rgs} shows the all such RGS of length 5
in lexicographic, co-lexicographic, subset-lexrev, and subset-lex
order.  The listing in subset-lex order does not seem to have
any apparently useful features, while the listing in subset-lexrev
order suggests the following method of generation.

We now assume a sentinel $a_n=+\infty$ at the end of the word.
\begin{alg}[Next-Catalan-Step-RGS]
Compute the successor in subset-lexrev order.
\end{alg}
\begin{enumerate}
\item Let $a_t$ be the first nonzero digit.
\item If $a_t < a_{t+1}$ and $a_t<t$, increment $a_t$ and return.
\item If $t\geq{}2$, set $a_{t-1}=1$ and return.
\item Remove all leading ones (setting them to zero),
 let $j$ be the position of the first digit $a_j\neq{}1$.
\item If the word is all-zero (that is, $j=n$), stop.
\item Set $a_j=a_j-1$ and set $a_{j-1}=1$.
\end{enumerate}

The implementation correctly handles all cases $n\geq{}0$.
{\codesize
\begin{listing}{1}
// FILE:  src/comb/catalan-step-rgs-subset-lexrev.h
class catalan_step_rgs_subset_lexrev
{
    ulong *a_;  // RGS
    ulong n2_;  // aux: min(n,2).
    ulong tr_;  // aux: track we are looking at
    ulong n_;   // length of RGS
\end{listing}
}%
The routine for the successor returns the
position of the rightmost change.
{\codesize
\begin{listing}{1}
    ulong next()
    {
        const ulong a0 = a_[tr_];

        if ( a0 < a_[tr_+1] )  // may read sentinel
        {
            if ( a0 < tr_ )  // can increment
            {
                a_[tr_] = a0 + 1;
                return tr_;
            }
        }

        if ( tr_ != 1 )  // can move left and increment (from 0 to 1)
        {
            --tr_;
            a_[tr_] = 1;
            return tr_;
        }

        // remove ones:
        ulong j = tr_;
        do  { a_[j] = 0; }  while ( a_[++j] == 1 );

        if ( j==n2_ )  return 0;  // current was last

        // decrement first value != 1:
        ulong aj = a_[j] - 1;
        a_[j] = aj;

        // move left and restore the 1:
        tr_ = j - 1;
        a_[tr_] = 1;
        return j;  // rightmost change at j==tr+1
    }
\end{listing}
}%

An update takes about 10.5 cycles,
whereas the updates for lexicographic and co-lexicographic
order respectively take 9 and 8 cycles.

%
%

%
\begin{figure}
{\jjbls
\begin{verbatim}
     0:  [ 1 1 1 1 1 1 1 1 1 1 ]      [ 1 1 1 1 1 1 1 1 1 1 ]
     1:  [ 1 1 1 1 1 1 1 1 2 . ]      [ 2 1 1 1 1 1 1 1 1 . ]
     2:  [ 1 1 1 1 1 1 1 3 . . ]      [ 2 2 1 1 1 1 1 1 . . ]
     3:  [ 1 1 1 1 1 1 2 2 . . ]      [ 3 1 1 1 1 1 1 1 . . ]
     4:  [ 1 1 1 1 1 1 4 . . . ]      [ 2 2 2 1 1 1 1 . . . ]
     5:  [ 1 1 1 1 1 2 3 . . . ]      [ 3 2 1 1 1 1 1 . . . ]
     6:  [ 1 1 1 1 2 2 2 . . . ]      [ 4 1 1 1 1 1 1 . . . ]
     7:  [ 1 1 1 1 1 5 . . . . ]      [ 2 2 2 2 1 1 . . . . ]
     8:  [ 1 1 1 1 2 4 . . . . ]      [ 3 2 2 1 1 1 . . . . ]
     9:  [ 1 1 1 1 3 3 . . . . ]      [ 3 3 1 1 1 1 . . . . ]
    10:  [ 1 1 1 2 2 3 . . . . ]      [ 4 2 1 1 1 1 . . . . ]
    11:  [ 1 1 2 2 2 2 . . . . ]      [ 5 1 1 1 1 1 . . . . ]
    12:  [ 1 1 1 1 6 . . . . . ]      [ 2 2 2 2 2 . . . . . ]
    13:  [ 1 1 1 2 5 . . . . . ]      [ 3 2 2 2 1 . . . . . ]
    14:  [ 1 1 1 3 4 . . . . . ]      [ 3 3 2 1 1 . . . . . ]
    15:  [ 1 1 2 2 4 . . . . . ]      [ 4 2 2 1 1 . . . . . ]
    16:  [ 1 1 2 3 3 . . . . . ]      [ 4 3 1 1 1 . . . . . ]
    17:  [ 1 2 2 2 3 . . . . . ]      [ 5 2 1 1 1 . . . . . ]
    18:  [ 2 2 2 2 2 . . . . . ]      [ 6 1 1 1 1 . . . . . ]
    19:  [ 1 1 1 7 . . . . . . ]      [ 3 3 2 2 . . . . . . ]
    20:  [ 1 1 2 6 . . . . . . ]      [ 4 2 2 2 . . . . . . ]
    21:  [ 1 1 3 5 . . . . . . ]      [ 3 3 3 1 . . . . . . ]
    22:  [ 1 2 2 5 . . . . . . ]      [ 4 3 2 1 . . . . . . ]
    23:  [ 1 1 4 4 . . . . . . ]      [ 5 2 2 1 . . . . . . ]
    24:  [ 1 2 3 4 . . . . . . ]      [ 4 4 1 1 . . . . . . ]
    25:  [ 2 2 2 4 . . . . . . ]      [ 5 3 1 1 . . . . . . ]
    26:  [ 1 3 3 3 . . . . . . ]      [ 6 2 1 1 . . . . . . ]
    27:  [ 2 2 3 3 . . . . . . ]      [ 7 1 1 1 . . . . . . ]
    28:  [ 1 1 8 . . . . . . . ]      [ 4 3 3 . . . . . . . ]
    29:  [ 1 2 7 . . . . . . . ]      [ 4 4 2 . . . . . . . ]
    30:  [ 1 3 6 . . . . . . . ]      [ 5 3 2 . . . . . . . ]
    31:  [ 2 2 6 . . . . . . . ]      [ 6 2 2 . . . . . . . ]
    32:  [ 1 4 5 . . . . . . . ]      [ 5 4 1 . . . . . . . ]
    33:  [ 2 3 5 . . . . . . . ]      [ 6 3 1 . . . . . . . ]
    34:  [ 2 4 4 . . . . . . . ]      [ 7 2 1 . . . . . . . ]
    35:  [ 3 3 4 . . . . . . . ]      [ 8 1 1 . . . . . . . ]
    36:  [ 1 9 . . . . . . . . ]      [ 5 5 . . . . . . . . ]
    37:  [ 2 8 . . . . . . . . ]      [ 6 4 . . . . . . . . ]
    38:  [ 3 7 . . . . . . . . ]      [ 7 3 . . . . . . . . ]
    39:  [ 4 6 . . . . . . . . ]      [ 8 2 . . . . . . . . ]
    40:  [ 5 5 . . . . . . . . ]      [ 9 1 . . . . . . . . ]
    41:  [10 . . . . . . . . . ]      [10 . . . . . . . . . ]
\end{verbatim}
}
\caption{\label{fig:part-asc-subset-lexrev}
The partitions of 10 in subset-lexrev order,
as weakly increasing lists of parts (left) and
as weakly decreasing lists of parts (right).
}
\end{figure}
%

We mention that the subset-lexrev order coincides
for certain objects with well-known orderings.
For example, for partitions as lists of parts
as shown in figure \ref{fig:part-asc-subset-lexrev},
both orderings are by falling length of partitions as
major order,
the minor order is lexicographic in case of weakly increasing parts
and reverse co-lexicographic for weakly decreasing parts.

Similar observations can be made (for example) for compositions into a fixed
number of parts, for both subset-lex and subset-lexrev order.

\section{SL-Gray order for binary words}

%
\begin{figure}
{\jjbls
\begin{verbatim}
    [ . . . . ]    [ . . . . ]    [ . . . . ]    [ . . . . ]
    [ . . . 1 ]    [ . . . 1 ]    [ . . . 1 ]    [ . . . 1 ]
    [ . . 1 . ]    [ . . 1 . ]    [ . . 1 . ] v  [ . . 1 1 ]
    [ . . 1 1 ]    [ . . 1 1 ]    [ . . 1 1 ] ^  [ . . 1 . ]
    [ . 1 . . ]    [ . 1 . . ] v  [ . 1 1 1 ] v  [ . 1 1 . ]
    [ . 1 . 1 ]    [ . 1 . 1 ]    [ . 1 1 . ] ^  [ . 1 1 1 ]
    [ . 1 1 . ]    [ . 1 1 . ]    [ . 1 . 1 ]    [ . 1 . 1 ]
    [ . 1 1 1 ]    [ . 1 1 1 ] ^  [ . 1 . . ]    [ . 1 . . ]
    [ 1 . . . ] v  [ 1 1 1 1 ] v  [ 1 1 . . ] v  [ 1 1 . . ]
    [ 1 . . 1 ]    [ 1 1 1 . ]    [ 1 1 . 1 ] ^  [ 1 1 . 1 ]
    [ 1 . 1 . ]    [ 1 1 . 1 ]    [ 1 1 1 . ] v  [ 1 1 1 1 ]
    [ 1 . 1 1 ]    [ 1 1 . . ] ^  [ 1 1 1 1 ] ^  [ 1 1 1 . ]
    [ 1 1 . . ]    [ 1 . 1 1 ]    [ 1 . 1 1 ] v  [ 1 . 1 . ]
    [ 1 1 . 1 ]    [ 1 . 1 . ]    [ 1 . 1 . ] ^  [ 1 . 1 1 ]
    [ 1 1 1 . ]    [ 1 . . 1 ]    [ 1 . . 1 ]    [ 1 . . 1 ]
    [ 1 1 1 1 ] ^  [ 1 . . . ]    [ 1 . . . ]    [ 1 . . . ]



  [ . . . . . ]    [ . . . . . ]    [ . . . . . ]    [ . . . . . ]
  [ 1 . . . . ]    [ 1 . . . . ]    [ 1 . . . . ]    [ 1 . . . . ]
  [ 1 1 . . . ]    [ 1 1 . . . ]    [ 1 1 . . . ]    [ 1 1 . . . ]
  [ 1 1 1 . . ]    [ 1 1 1 . . ]    [ 1 1 1 . . ]    [ 1 1 1 . . ]
  [ 1 1 1 1 . ]    [ 1 1 1 1 . ]    [ 1 1 1 1 . ]    [ 1 1 1 1 . ]
  [ 1 1 1 1 1 ]    [ 1 1 1 1 1 ]    [ 1 1 1 1 1 ]    [ 1 1 1 1 1 ]
  [ 1 1 1 . 1 ]    [ 1 1 1 . 1 ]    [ 1 1 1 . 1 ]    [ 1 1 1 . 1 ]
  [ 1 1 . 1 . ]    [ 1 1 . 1 . ]    [ 1 1 . 1 . ] v  [ 1 1 . . 1 ]
  [ 1 1 . 1 1 ]    [ 1 1 . 1 1 ]    [ 1 1 . 1 1 ]    [ 1 1 . 1 1 ]
  [ 1 1 . . 1 ]    [ 1 1 . . 1 ]    [ 1 1 . . 1 ] ^  [ 1 1 . 1 . ]
  [ 1 . 1 . . ]    [ 1 . 1 . . ] v  [ 1 . . . 1 ] v  [ 1 . . 1 . ]
  [ 1 . 1 1 . ]    [ 1 . 1 1 . ]    [ 1 . . 1 1 ]    [ 1 . . 1 1 ]
  [ 1 . 1 1 1 ]    [ 1 . 1 1 1 ]    [ 1 . . 1 . ] ^  [ 1 . . . 1 ]
  [ 1 . 1 . 1 ]    [ 1 . 1 . 1 ]    [ 1 . 1 . 1 ]    [ 1 . 1 . 1 ]
  [ 1 . . 1 . ]    [ 1 . . 1 . ]    [ 1 . 1 1 1 ]    [ 1 . 1 1 1 ]
  [ 1 . . 1 1 ]    [ 1 . . 1 1 ]    [ 1 . 1 1 . ]    [ 1 . 1 1 . ]
  [ 1 . . . 1 ]    [ 1 . . . 1 ] ^  [ 1 . 1 . . ]    [ 1 . 1 . . ]
  [ . 1 . . . ] v  [ . . . . 1 ] v  [ . . 1 . . ]    [ . . 1 . . ]
  [ . 1 1 . . ]    [ . . . 1 1 ]    [ . . 1 1 . ]    [ . . 1 1 . ]
  [ . 1 1 1 . ]    [ . . . 1 . ]    [ . . 1 1 1 ]    [ . . 1 1 1 ]
  [ . 1 1 1 1 ]    [ . . 1 . 1 ]    [ . . 1 . 1 ]    [ . . 1 . 1 ]
  [ . 1 1 . 1 ]    [ . . 1 1 1 ]    [ . . . 1 . ] v  [ . . . . 1 ]
  [ . 1 . 1 . ]    [ . . 1 1 . ]    [ . . . 1 1 ]    [ . . . 1 1 ]
  [ . 1 . 1 1 ]    [ . . 1 . . ] ^  [ . . . . 1 ] ^  [ . . . 1 . ]
  [ . 1 . . 1 ]    [ . 1 . . 1 ]    [ . 1 . . 1 ] v  [ . 1 . 1 . ]
  [ . . 1 . . ]    [ . 1 . 1 1 ]    [ . 1 . 1 1 ]    [ . 1 . 1 1 ]
  [ . . 1 1 . ]    [ . 1 . 1 . ]    [ . 1 . 1 . ] ^  [ . 1 . . 1 ]
  [ . . 1 1 1 ]    [ . 1 1 . 1 ]    [ . 1 1 . 1 ]    [ . 1 1 . 1 ]
  [ . . 1 . 1 ]    [ . 1 1 1 1 ]    [ . 1 1 1 1 ]    [ . 1 1 1 1 ]
  [ . . . 1 . ]    [ . 1 1 1 . ]    [ . 1 1 1 . ]    [ . 1 1 1 . ]
  [ . . . 1 1 ]    [ . 1 1 . . ]    [ . 1 1 . . ]    [ . 1 1 . . ]
  [ . . . . 1 ] ^  [ . 1 . . . ]    [ . 1 . . . ]    [ . 1 . . . ]
\end{verbatim}
}
\caption{\label{fig:gray-rec}
Construction of Gray codes by reversing sublists:
for the binary reflected Gray code (top)
and for the binary SL-Gray order (bottom).
The symbols \texttt{v} and \texttt{\^{}} respectively mark
begin and end of the reversed sublists.
}
\end{figure}

A well-known construction for the binary reflected Gray code proceeds
by successively reversing ranges having identical prefixes that
end in a one, using increasingly longer prefixes,
see top of figure \ref{fig:gray-rec}.

The same construction, now using prefixes ending in a zero,
gives a Gray code for subset-lex order, see
bottom of figure \ref{fig:gray-rec}.
We will call this ordering the \jjterm{SL-Gray order}.
%

For the computation of the successor we keep
a variable $t$ for the current track and a variable
$d$ indicating the direction in which the track
will be moved if necessary.
Initially $t=0$ and $d=+1$.
\begin{samepage}
\begin{alg}[Next-SL-Gray]
Compute the successor in SL-Gray order.
\end{alg}
\begin{enumerate}
\item If $d=+1$, do the following (try to append trailing ones):
\begin{enumerate}
\item If $a_t=0$, set $a_t=1$, $t=t+1$, and return.
\item Otherwise, set $d=-1$ (change direction),
  $t=n-1$ (move to rightmost track),
  $j=n-2$, $a_j=1-a_j$, and return.
\end{enumerate}
\item Otherwise ($d=-1$), do the following (try to remove trailing ones):
\begin{enumerate}
\item If $a_{t-1}=1$, set $a_t=0$, $t=t-1$, and return.
\item Otherwise, set $d=+1$ (change direction),
  $j=t-2$,  $a_j=1-a_j$,  $t=t+1$ (move right), and return.
\end{enumerate}
\end{enumerate}
\end{samepage}

The algorithm is loopless.
In the implementation, the variables
\texttt{tr} and \texttt{dt} respectively
correspond to $t$ and $d$ in the algorithm.
Two sentinels $a_{-1}=a_n=+1$ are used.
{\codesize
\begin{listing}{1}
// FILE:  src/comb/binary-sl-gray.h
class binary_sl_gray
{
    ulong n_;   // number of digits
    ulong tr_;  // aux: current track (0 <= tr <= n)
    ulong dt_;  // aux: direction in which track tries to move
    ulong *a_;  // digits

    void first()
    {
        for (ulong k=0; k<n_; ++k)  a_[k] = 0;
        tr_ = 0;
        dt_ = +1;
    }

    explicit binary_sl_gray(ulong n)
    {
        n_ = n;
        a_ = new ulong[n_+2];  // sentinels at both ends

        a_[n_+1] = +1;  // != 0
        a_[0] = +1;

        ++a_;  // nota bene

        first();
    }

    void first()
    {
        for (ulong k=0; k<n_; ++k)  a_[k] = 0;
        tr_ = 0;
        dt_ = +1;

        j_ = ( n_>=2 ? 1 : 0);  // wrt. last word
        dm_ = -1;  // wrt. last word
    }

\end{listing}
}%
The routines for successor and predecessor handle
all cases $n\geq{}0$ correctly.
{\codesize
\begin{listing}{1}
    bool next()
    {
        if ( dt_ == +1 )  // try to append trailing ones
        {
            if ( a_[tr_] == 0 )  // can append  // may read sentinel a[n]
            {
                a_[tr_] = 1;
                ++tr_;
            }
            else
            {
                dt_ = -1UL;  // change direction
                tr_ = n_ - 1;
                ulong j_ = tr_ - 1;
                // Is current last (only for n_ <= 1)?
                if ( j_ > n_ )  return false;
                a_[j_] = 1 - a_[j_];
            }
        }
        else  // dt_ == -1  // try to remove trailing ones
        {
            if ( a_[tr_-1] != 0 )   // can remove  // tr - 1 >= 0
            {
                a_[tr_] = 0;
                --tr_;
            }
            else
            {
                dt_ = +1;  // change direction
                ulong j_ = tr_ - 2;
                if ( (long)j_ < 0 )  return false;
                a_[j_] = 1 - a_[j_];
                ++tr_;
            }
        }

        return true;
    }
\end{listing}
}%
The routine to compute the predecessor is obtained
by (essentially) negating the direction the
track tries to move (variable \texttt{dt}):
{\codesize
\begin{listing}{1}
    bool prev()
    {
        if ( dt_ != +1 )  // dt==-1  // try to append trailing ones
        {
            if ( a_[tr_+1] == 0 )  // can append  // may read sentinel a[n]
            {
                a_[tr_+1] = 1;
                ++tr_;
            }
            else
            {
                dt_ = +1;  // change direction
                tr_ = n_ - 1;
                ulong j_ = tr_ - 1;
                a_[j_] = 1 - a_[j_];
                ++tr_;
            }
        }
        else  // dt_ == +1  // try to remove trailing ones
        {
            if ( tr_ == 0 )
            {
                if ( a_[0]==0 )  return false;  // (only for n_ <= 1)
                a_[0] = 0;
                return true;
            }

            // tr - 1 >= -1 (can read low sentinel)
            if ( a_[tr_-2] != 0 )  // can remove
            {
                a_[tr_-1] = 0;
                --tr_;
            }
            else
            {
                dt_ = -1UL;  // change direction
                ulong j_ = tr_ - 3;
                a_[j_] = 1 - a_[j_];
                --tr_;
            }
        }

        return true;
    }
\end{listing}
}%
One update in either direction takes about 7.5 cycles.
An implementation for the generation of the
SL-Gray order in a binary word is
given in \fileref{src/bits/bit-sl-gray.h}.

The ranking algorithm is easily obtained by
observing that if the highest bit is not set
(and the word is nonzero),
the order for the remaining word is reflected,
as shown in figure \ref{fig:gray-rec}.
\begin{alg}[Binary-SL-Gray-Rank]
Recursive routine $F(k,\,n)$ for the computation of the rank of
 an $n$-bit word $k$ in binary SL-Gray order.
\end{alg}
\begin{enumerate}
\item If $k=0$, return 0.
\item Set $w=k$ and unset the highest bit in $w$.
\item If the highest bit of $k$ (at position $n-1$) is set,
 return $1 + F(w,\, n-1)$.
\item Otherwise return $2^n - F(w,\, n-1)$.
\end{enumerate}
%
In the following implementation
the variable \texttt{ldn} must give
the number of bits in the SL-Gray code.
{\codesize
\begin{listing}{1}
// FILE:  src/bits/bin-to-sl-gray.h
ulong sl_gray_to_bin(ulong k, ulong ldn)
{
    if ( k==0 )  return 0;
    ulong b = 1UL << (ldn-1);  // mask for bit at end
    ulong h = k & b;  // bit at end
    k ^= h;
    ulong z = sl_gray_to_bin( k, ldn-1 );  // recursion
    if ( h==0 )  return (b<<1) - z;
    else         return 1 + z;
}
\end{listing}
}%

A routine for the conversion of a binary
word into the corresponding word in SL-Gray order
(unranking algorithm) is
{\codesize
\begin{listing}{1}
ulong bin_to_sl_gray(ulong k, ulong ldn)
{
    if ( ldn==0 )  return 0;

    ulong b = 1UL << (ldn-1);  // highest bit
    ulong m = (b<<1) - 1;  // mask for reversing direction
    ulong z = b;  // Gray code
    k -= 1;  // move all-zero word to begin
    while ( b != 0 )
    {
        const ulong h = k & b;  // bit under consideration
        z ^= h;  // with one, switch bit in Gray code
        if ( !h )  k ^= m;  // reverse direction with zero
        k += 1;  // SL-Gray
        b >>= 1;  // next lower bit
        m >>= 1;  // next smaller mask
    }

    return z;
}
\end{listing}
}%
%
%

\begin{figure}
{\jjbls
\begin{verbatim}


        0121030121014101210301210125210121030121014101210301210123

  ......1...1.1...1...1...1.1...1.....1...1.1...1...1...1.1...1...
  .....1.1.1...1.1.1.1.1.1...1.1.1...1.1.1...1.1.1.1.1.1...1.1.1..
  ....1...1.....1.......1.....1...1.1...1.....1.......1.....1...1.
  ...1.......1.............1...............1.............1.......1
  ..1...............1.............................1...............
  .1...............................1..............................

  ......1111..1111....1111..1111......1111..1111....1111..1111....
  .....11..1111..11..11..1111..11....11..1111..11..11..1111..11...
  ....1111......11111111......1111..1111......11111111......1111..
  ...11111111..............1111111111111111..............11111111.
  ..1111111111111111..............................1111111111111111
  .11111111111111111111111111111111...............................

\end{verbatim}
}
\caption{\label{fig:bin-sl-gray-delta}
The delta sequence (top) of the binary SL-Gray order,
starting after the initial slope
and indexing positions from the end of the words.
}
\end{figure}
%

The delta sequence (sequence of positions of changes),
starting after the initial slope and indexing positions from the end,
is shown in figure \ref{fig:bin-sl-gray-delta},
this is sequence \jjseqref{A217262} in \cite{oeis}.
It can be obtained from the ruler function (sequence \jjseqref{A007814})
by replacing $0$ by $w=01210$,
 $1$ by $3$, $2$ by $141$, $3$ by $12521$, $4$ by $1236321$, \ldots,
 $n$ by $123\ldots{}(n-1)(n+2)(n-1)\ldots{}321$:
{\jjbls
\begin{verbatim}
0 1 0  2  0 1 0   3   0 1 0  2  0 1 0    4    0 ...  (ruler function)
w 3 w 141 w 3 w 12521 w 3 w 141 w 3 w 1236321 w ... =

01210 3 01210 141 01210 3 01210 12521 01210 3 01210 141 01210 3 01210 ...
\end{verbatim}
}%
%
It is easy to see that (for $n\geq{}4$) the
SL-Gray order for the $2^n$ words of length $n$
has $2^{n-2}-2$ successive transitions that are 3-close
(have a distance of 3), the rest are 1-close
(adjacent changes).

A Gray code where the maximal distance between successive
transitions is $k$ is called $k$-close.
In \cite[p.400]{fxtbook} it has been observed
that 1-close Gray codes exists for $n\leq{}6$
but not for $n=7$ (and it appears unlikely that
for any $n\geq{}8$ such a Gray code exists).
Does a 2-close Gray code exist?
We remark that the SL-Gray order is 2-close
if the radices for all digits are even and $\geq{}4$.


\section{SL-Gray order for mixed radix words}

%
\begin{figure}
\vspace*{-8mm}
{\jjbls
\begin{verbatim}
   0:  [ . . . . ]    [ + + + + ]    {  }
   1:  [ 1 . . . ]    [ + + + + ]    { 0 }
   2:  [ 1 1 . . ]    [ - + + + ]    { 0, 1 }
   3:  [ 1 2 . . ]    [ - + + + ]    { 0, 1, 1 }
   4:  [ 1 2 1 . ]    [ - - + + ]    { 0, 1, 1, 2 }
   5:  [ 1 2 2 . ]    [ - - + + ]    { 0, 1, 1, 2, 2 }
   6:  [ 1 2 2 1 ]    [ - - - + ]    { 0, 1, 1, 2, 2, 3 }
   7:  [ 1 2 2 2 ]    [ - - - + ]    { 0, 1, 1, 2, 2, 3, 3 }
   8:  [ 1 2 2 3 ]    [ - - - + ]    { 0, 1, 1, 2, 2, 3, 3, 3 }
   9:  [ 1 2 1 3 ]    [ - - - - ]    { 0, 1, 1, 2, 3, 3, 3 }
  10:  [ 1 2 1 2 ]    [ - - - - ]    { 0, 1, 1, 2, 3, 3 }
  11:  [ 1 2 1 1 ]    [ - - - - ]    { 0, 1, 1, 2, 3 }
  12:  [ 1 2 . 1 ]    [ - - - + ]    { 0, 1, 1, 3 }
  13:  [ 1 2 . 2 ]    [ - - - + ]    { 0, 1, 1, 3, 3 }
  14:  [ 1 2 . 3 ]    [ - - - + ]    { 0, 1, 1, 3, 3, 3 }
  15:  [ 1 1 . 3 ]    [ - - + - ]    { 0, 1, 3, 3, 3 }
  16:  [ 1 1 . 2 ]    [ - - + - ]    { 0, 1, 3, 3 }
  17:  [ 1 1 . 1 ]    [ - - + - ]    { 0, 1, 3 }
  18:  [ 1 1 1 1 ]    [ - - + + ]    { 0, 1, 2, 3 }
  19:  [ 1 1 1 2 ]    [ - - + + ]    { 0, 1, 2, 3, 3 }
  20:  [ 1 1 1 3 ]    [ - - + + ]    { 0, 1, 2, 3, 3, 3 }
  21:  [ 1 1 2 3 ]    [ - - + - ]    { 0, 1, 2, 2, 3, 3, 3 }
  22:  [ 1 1 2 2 ]    [ - - + - ]    { 0, 1, 2, 2, 3, 3 }
  23:  [ 1 1 2 1 ]    [ - - + - ]    { 0, 1, 2, 2, 3 }
  24:  [ 1 1 2 . ]    [ - - - + ]    { 0, 1, 2, 2 }
  25:  [ 1 1 1 . ]    [ - - - + ]    { 0, 1, 2 }
  26:  [ 1 . 1 . ]    [ - - + + ]    { 0, 2 }
  27:  [ 1 . 2 . ]    [ - - + + ]    { 0, 2, 2 }
  28:  [ 1 . 2 1 ]    [ - - - + ]    { 0, 2, 2, 3 }
  29:  [ 1 . 2 2 ]    [ - - - + ]    { 0, 2, 2, 3, 3 }
  30:  [ 1 . 2 3 ]    [ - - - + ]    { 0, 2, 2, 3, 3, 3 }
  31:  [ 1 . 1 3 ]    [ - - - - ]    { 0, 2, 3, 3, 3 }
  32:  [ 1 . 1 2 ]    [ - - - - ]    { 0, 2, 3, 3 }
  33:  [ 1 . 1 1 ]    [ - - - - ]    { 0, 2, 3 }
  34:  [ 1 . . 1 ]    [ - - - + ]    { 0, 3 }
  35:  [ 1 . . 2 ]    [ - - - + ]    { 0, 3, 3 }
  36:  [ 1 . . 3 ]    [ - - - + ]    { 0, 3, 3, 3 }
  37:  [ . . . 3 ]    [ - + + - ]    { 3, 3, 3 }
  38:  [ . . . 2 ]    [ - + + - ]    { 3, 3 }
  39:  [ . . . 1 ]    [ - + + - ]    { 3 }
  40:  [ . . 1 1 ]    [ - + + + ]    { 2, 3 }
  41:  [ . . 1 2 ]    [ - + + + ]    { 2, 3, 3 }
  42:  [ . . 1 3 ]    [ - + + + ]    { 2, 3, 3, 3 }
  43:  [ . . 2 3 ]    [ - + + - ]    { 2, 2, 3, 3, 3 }
  44:  [ . . 2 2 ]    [ - + + - ]    { 2, 2, 3, 3 }
  45:  [ . . 2 1 ]    [ - + + - ]    { 2, 2, 3 }
  46:  [ . . 2 . ]    [ - + - + ]    { 2, 2 }
  47:  [ . . 1 . ]    [ - + - + ]    { 2 }
  48:  [ . 1 1 . ]    [ - + + + ]    { 1, 2 }
  49:  [ . 1 2 . ]    [ - + + + ]    { 1, 2, 2 }
  50:  [ . 1 2 1 ]    [ - + - + ]    { 1, 2, 2, 3 }
  51:  [ . 1 2 2 ]    [ - + - + ]    { 1, 2, 2, 3, 3 }
  52:  [ . 1 2 3 ]    [ - + - + ]    { 1, 2, 2, 3, 3, 3 }
  53:  [ . 1 1 3 ]    [ - + - - ]    { 1, 2, 3, 3, 3 }
  54:  [ . 1 1 2 ]    [ - + - - ]    { 1, 2, 3, 3 }
  55:  [ . 1 1 1 ]    [ - + - - ]    { 1, 2, 3 }
  56:  [ . 1 . 1 ]    [ - + - + ]    { 1, 3 }
  57:  [ . 1 . 2 ]    [ - + - + ]    { 1, 3, 3 }
  58:  [ . 1 . 3 ]    [ - + - + ]    { 1, 3, 3, 3 }
  59:  [ . 2 . 3 ]    [ - + + - ]    { 1, 1, 3, 3, 3 }
  60:  [ . 2 . 2 ]    [ - + + - ]    { 1, 1, 3, 3 }
  61:  [ . 2 . 1 ]    [ - + + - ]    { 1, 1, 3 }
  62:  [ . 2 1 1 ]    [ - + + + ]    { 1, 1, 2, 3 }
  63:  [ . 2 1 2 ]    [ - + + + ]    { 1, 1, 2, 3, 3 }
  64:  [ . 2 1 3 ]    [ - + + + ]    { 1, 1, 2, 3, 3, 3 }
  65:  [ . 2 2 3 ]    [ - + + - ]    { 1, 1, 2, 2, 3, 3, 3 }
  66:  [ . 2 2 2 ]    [ - + + - ]    { 1, 1, 2, 2, 3, 3 }
  67:  [ . 2 2 1 ]    [ - + + - ]    { 1, 1, 2, 2, 3 }
  68:  [ . 2 2 . ]    [ - + - + ]    { 1, 1, 2, 2 }
  69:  [ . 2 1 . ]    [ - + - + ]    { 1, 1, 2 }
  70:  [ . 2 . . ]    [ - - + + ]    { 1, 1 }
  71:  [ . 1 . . ]    [ - - + + ]    { 1 }
\end{verbatim}
}
\caption{\label{fig:mixedradix-sl-gray}
Mixed radix words (left) and corresponding subsets (right)
of the multiset $\{0^1, 1^2, 2^2, 3^3 \}$ in SL-Gray order,
and the array of directions used in the computation (middle).
}
\end{figure}
%

The generalization of the SL-Gray order to mixed radix words is obtained by
successive reversion of the sublists with identical prefix for all prefixes
ending in an even digit.

For the computation of the successor we keep
a variable $t$ for the current track and
variables $d_j=\pm{}1$ (an array) for the direction
the digit $a_j$ is currently moving in
(whether its is incremented or decremented).
We further use $m_j$ for the maximal value
the digit $a_j$ can have (the \eqq{nines}).
Initially all digits $a_k$ are zero,
all directions $d_k$ are $+1$, and $t=0$.
\begin{alg}[Next-SL-Gray]
Compute the successor SL-Gray order.
\end{alg}
\begin{enumerate}
\item Set $j=t$ and $b = a_j + d_j$.
\item If $b\neq{}0$ and $b\leq{}m_j$, set $a_j=b$ and return (easy case).
\item Set $d_j=-d_j$ (change direction for digit $j$).
\item If $d_j = +1$ and $a_{j+1}=0$,
 set $a_{j+1} = +1$, $t = j+1$ (move track right), and return.
\item If $d_j=-1$ and $a_{j-1} = m_{j-1}$,
 set $a_j=0$, $d_{j-1} = -1$, $t=j-1$ (move track left), and return.
\item Find the position $p$ of the nearest digit $a_k$ to the left
 such that $a_k+d_k$ is in the range $0,1,\ldots,m_k$ (a valid digit).
 In the process, change $d_k$ for all $k$ where $p < k < j$.
\item Set $a_p = a_p + d_p$ (change digit, keep track).
\end{enumerate}
%
The details for termination have been omitted,
this is handled with the help of sentinels
in the implementation.
%
%
{\codesize
\begin{listing}{1}
// FILE: src/comb/mixedradix-sl-gray.h
class mixedradix_sl_gray
{
    ulong n_;    // number of digits
                 // (n kinds of elements in multiset, n>=1)
    ulong tr_;   // aux: current track
    ulong *a_;   // digits of mixed radix number
                 // (multiplicity of kind k in subset).
    ulong *d_;   // directions (either +1 or -1)
    ulong *m1_;  // nines (radix minus one) for each digit
                 // (multiplicity of kind k in set).
\end{listing}
}%
Sentinels are used in all arrays at index $-1$ and at index $n$
to handle termination.
{\codesize
\begin{listing}{1}
    mixedradix_sl_gray(ulong n, ulong mm, const ulong *m=0)
    // Must have n>=1
    {
        n_ = n;
        a_ = new ulong[n_+2];  // all with sentinels at both ends
        d_ = new ulong[n_+2];
        m1_ = new ulong[n_+2];

        // sentinels on the right:
        a_[n_+1] = +1;   // != 0
        m1_[n_+1] = +1;  // same as a_[n+1]
        d_[n_+1] = +1;   // positive

        // sentinels on the left:
        a_[0] = +2;   // >= +2
        m1_[0] = +2;  //  same as a_[0]
        d_[0] = 0;    // zero

        ++a_;  ++d_;  ++m1_;  // nota bene

        mixedradix_init(n_, mm, m, m1_);
        first();
    }

    void first()
    {
        for (ulong k=0; k<n_; ++k)  a_[k] = 0;
        for (ulong k=0; k<n_; ++k)  d_[k] = +1;
        tr_ = 0;
    }
\end{listing}
}%
The method for the successor handles all cases $n\geq{}0$ correctly.
The variable \texttt{a1} in the implementation corresponds to $b$ in the algorithm.
{\codesize
\begin{listing}{1}
    bool next()
    {
        ulong j = tr_;
        const ulong dj = d_[j];
        const ulong a1 = a_[j] + dj;  // a[j] +- 1

        if ( (a1 != 0) && (a1 <= m1_[j]) )  // easy case
        {
            a_[j] = a1;
            return true;
        }

        d_[j] = -dj;  // change direction

        if ( dj == +1 )  // so a_[j] == m1_[j] == nine
        {
            // Try to move track right with a[j] == nine:
            const ulong j1 = j + 1;
            if ( a_[j1] == 0 )  // can read high sentinel
            {
                a_[j1] = +1;
                tr_ = j1;
                return true;
            }
        }
        else  // here dj == -1, so a_[j] == 1
        {
            if ( (long)j <= 0 )  return false;  // current is last

            // Try to move track left with a[j] == 1:
            const ulong j1 = j - 1;
            if ( a_[j1] == m1_[j1] )  // can read low sentinel when n_ == 1
            {
                a_[j] = 0;
                d_[j1] = -1UL;
                tr_ = j1;
                return true;
            }
        }

        // find first changeable track to the left:
        --j;
        while ( a_[j] + d_[j] > m1_[j] )  // may read low sentinels
        {
            d_[j] = -d_[j];  // change direction
            --j;
        }

        if ( (long)j < 0 )  return false;  // current is last

        // Change digit left but keep track:
        a_[j] += d_[j];

        return true;
    }
\end{listing}
}%
%


%
%
%

\section{A Gray code for compositions}

%
\begin{figure}
{\jjbls
\begin{verbatim}

   0:  111111  [ 7 ]
   1:  11111.  [ 6 1 ]
   2:  1111..  [ 5 1 1 ]                   0:  111.....  [ 4 1 1 1 1 1 ]
   3:  1111.1  [ 5 2 ]                     1:  11..1...  [ 3 1 2 1 1 1 ]
   4:  111.11  [ 4 3 ]                     2:  11....1.  [ 3 1 1 1 2 1 ]
   5:  111.1.  [ 4 2 1 ]                   3:  11.....1  [ 3 1 1 1 1 2 ]
   6:  111...  [ 4 1 1 1 ]                 4:  11...1..  [ 3 1 1 2 1 1 ]
   7:  111..1  [ 4 1 2 ]                   5:  11.1....  [ 3 2 1 1 1 1 ]
   8:  11..11  [ 3 1 3 ]                   6:  1.11....  [ 2 3 1 1 1 1 ]
   9:  11..1.  [ 3 1 2 1 ]                 7:  1.1.1...  [ 2 2 2 1 1 1 ]
  10:  11....  [ 3 1 1 1 1 ]               8:  1.1...1.  [ 2 2 1 1 2 1 ]
  11:  11...1  [ 3 1 1 2 ]                 9:  1.1....1  [ 2 2 1 1 1 2 ]
  12:  11.1..  [ 3 2 1 1 ]                10:  1.1..1..  [ 2 2 1 2 1 1 ]
  13:  11.1.1  [ 3 2 2 ]                  11:  1...11..  [ 2 1 1 3 1 1 ]
  14:  11.11.  [ 3 3 1 ]                  12:  1...1.1.  [ 2 1 1 2 2 1 ]
  15:  11.111  [ 3 4 ]                    13:  1...1..1  [ 2 1 1 2 1 2 ]
  16:  1.1111  [ 2 5 ]                    14:  1.....11  [ 2 1 1 1 1 3 ]
  17:  1.111.  [ 2 4 1 ]                  15:  1....1.1  [ 2 1 1 1 2 2 ]
  18:  1.11..  [ 2 3 1 1 ]                16:  1....11.  [ 2 1 1 1 3 1 ]
  19:  1.11.1  [ 2 3 2 ]                  17:  1..1..1.  [ 2 1 2 1 2 1 ]
  20:  1.1.11  [ 2 2 3 ]                  18:  1..1...1  [ 2 1 2 1 1 2 ]
  21:  1.1.1.  [ 2 2 2 1 ]                19:  1..1.1..  [ 2 1 2 2 1 1 ]
  22:  1.1...  [ 2 2 1 1 1 ]              20:  1..11...  [ 2 1 3 1 1 1 ]
  23:  1.1..1  [ 2 2 1 2 ]                21:  ..111...  [ 1 1 4 1 1 1 ]
  24:  1...11  [ 2 1 1 3 ]                22:  ..11..1.  [ 1 1 3 1 2 1 ]
  25:  1...1.  [ 2 1 1 2 1 ]              23:  ..11...1  [ 1 1 3 1 1 2 ]
  26:  1.....  [ 2 1 1 1 1 1 ]            24:  ..11.1..  [ 1 1 3 2 1 1 ]
  27:  1....1  [ 2 1 1 1 2 ]              25:  ..1.11..  [ 1 1 2 3 1 1 ]
  28:  1..1..  [ 2 1 2 1 1 ]              26:  ..1.1.1.  [ 1 1 2 2 2 1 ]
  29:  1..1.1  [ 2 1 2 2 ]                27:  ..1.1..1  [ 1 1 2 2 1 2 ]
  30:  1..11.  [ 2 1 3 1 ]                28:  ..1...11  [ 1 1 2 1 1 3 ]
  31:  1..111  [ 2 1 4 ]                  29:  ..1..1.1  [ 1 1 2 1 2 2 ]
  32:  ..1111  [ 1 1 5 ]                  30:  ..1..11.  [ 1 1 2 1 3 1 ]
  33:  ..111.  [ 1 1 4 1 ]                31:  ....111.  [ 1 1 1 1 4 1 ]
  34:  ..11..  [ 1 1 3 1 1 ]              32:  ....11.1  [ 1 1 1 1 3 2 ]
  35:  ..11.1  [ 1 1 3 2 ]                33:  ....1.11  [ 1 1 1 1 2 3 ]
  36:  ..1.11  [ 1 1 2 3 ]                34:  .....111  [ 1 1 1 1 1 4 ]
  37:  ..1.1.  [ 1 1 2 2 1 ]              35:  ...1..11  [ 1 1 1 2 1 3 ]
  38:  ..1...  [ 1 1 2 1 1 1 ]            36:  ...1.1.1  [ 1 1 1 2 2 2 ]
  39:  ..1..1  [ 1 1 2 1 2 ]              37:  ...1.11.  [ 1 1 1 2 3 1 ]
  40:  ....11  [ 1 1 1 1 3 ]              38:  ...11.1.  [ 1 1 1 3 2 1 ]
  41:  ....1.  [ 1 1 1 1 2 1 ]            39:  ...11..1  [ 1 1 1 3 1 2 ]
  42:  ......  [ 1 1 1 1 1 1 1 ]          40:  ...111..  [ 1 1 1 4 1 1 ]
  43:  .....1  [ 1 1 1 1 1 2 ]            41:  .1..11..  [ 1 2 1 3 1 1 ]
  44:  ...1..  [ 1 1 1 2 1 1 ]            42:  .1..1.1.  [ 1 2 1 2 2 1 ]
  45:  ...1.1  [ 1 1 1 2 2 ]              43:  .1..1..1  [ 1 2 1 2 1 2 ]
  46:  ...11.  [ 1 1 1 3 1 ]              44:  .1....11  [ 1 2 1 1 1 3 ]
  47:  ...111  [ 1 1 1 4 ]                45:  .1...1.1  [ 1 2 1 1 2 2 ]
  48:  .1..11  [ 1 2 1 3 ]                46:  .1...11.  [ 1 2 1 1 3 1 ]
  49:  .1..1.  [ 1 2 1 2 1 ]              47:  .1.1..1.  [ 1 2 2 1 2 1 ]
  50:  .1....  [ 1 2 1 1 1 1 ]            48:  .1.1...1  [ 1 2 2 1 1 2 ]
  51:  .1...1  [ 1 2 1 1 2 ]              49:  .1.1.1..  [ 1 2 2 2 1 1 ]
  52:  .1.1..  [ 1 2 2 1 1 ]              50:  .1.11...  [ 1 2 3 1 1 1 ]
  53:  .1.1.1  [ 1 2 2 2 ]                51:  .11.1...  [ 1 3 2 1 1 1 ]
  54:  .1.11.  [ 1 2 3 1 ]                52:  .11...1.  [ 1 3 1 1 2 1 ]
  55:  .1.111  [ 1 2 4 ]                  53:  .11....1  [ 1 3 1 1 1 2 ]
  56:  .11.11  [ 1 3 3 ]                  54:  .11..1..  [ 1 3 1 2 1 1 ]
  57:  .11.1.  [ 1 3 2 1 ]                55:  .111....  [ 1 4 1 1 1 1 ]
  58:  .11...  [ 1 3 1 1 1 ]
  59:  .11..1  [ 1 3 1 2 ]
  60:  .111..  [ 1 4 1 1 ]
  61:  .111.1  [ 1 4 2 ]
  62:  .1111.  [ 1 5 1 ]
  63:  .11111  [ 1 6 ]
\end{verbatim}
}
\caption{\label{fig:comp-gray}
The 64 compositions of 7 together with their binary encodings
as in figure \ref{fig:comp-nz} (left)
and the 56 compositions of 9 into exactly 6 parts (right).
}
\end{figure}
%

We now describe a Gray code for compositions where at most one unit is
moved in each step, all moves are 1-close or 2-close (these always
cross a part 1), and all moves are at the end of the current composition.

The compositions of 7 in this ordering are shown in
figure \ref{fig:comp-gray} (left), it can be obtained
from the lexicographic order by successively reversing the sublists
with identical prefixes that end in an odd part.
To keep matters simple for the iterative algorithm, we arrange the compositions
for $n$ odd such that the first part is only decreasing (starting with the
composition $[n]$) and otherwise such that the first part is increasing
(starting with the composition $[1,n-1]$).  The compositions of 6 in this
ordering are obtained by dropping the first part 1 in all compositions of 7 in
the second half of the list, see figure \ref{fig:comp-gray}.

A recursive routine can be obtained by switching
between the recursive functions for lexicographic
and reversed lexicographic order whenever
an odd part has been written.
In the following a global array \texttt{a[]} is used.
{\codesize
\begin{listing}{1}
// FILE:  src/comb/composition-nz-gray-rec-demo.cc
void F(ulong n, ulong m)
{
    if ( n==0 )  { visit(m);  return; }

    for (ulong f=n; f!=0; --f)  // first part decreasing
    {
        a[m] = f;
        if ( 0 == (f & 1) )  F(n-f, m+1);
        else                 B(n-f, m+1);
    }
}

void B(ulong n, ulong m)
{
    if ( n==0 )  { visit(m);  return; }

    for (ulong f=1; f<=n; ++f)  // first part increasing
    {
        a[m] = f;
        if ( 0 == (f & 1) )  B(n-f, m+1);
        else                 F(n-f, m+1);
    }
}
\end{listing}
}%
The initial call is \texttt{F(n, 0)} for $n$ odd
and \texttt{G(n,0)} for $n$ even.

For the iterative algorithm we use a function
$D(x)$ that shall return $+1$ if $x$ is odd and otherwise $-1$.
We will refer to the last three parts as respectively $x$, $y$, and $z$.
\begin{alg}[Next-Comp-Gray]
Compute the successor of a composition.
\end{alg}
\begin{enumerate}
\item If $z=n-1$ and $n$ is odd, or $z=n$ and $n$ is even, stop.
%
\item If $z=1$, do the following:
\begin{enumerate}
\item If $D(y)=+1$, set $y=y+1$ and drop $z$; return.
\item Otherwise ($D(y)=-1$), set $y=y-1$ and append a part $1$; return.
\end{enumerate}
%
\item Otherwise ($z>1$) do the following:
\begin{enumerate}
\item If $z$ is odd, set $z=z-1$ and append a part $1$; return.
\item If $y>1$, set $y=y-D(y)$ and $z=z+D(y)$; return.
\item If $x>1$, set $x=x+D(x)$ and $z=z-D(x)$; return.
\item Otherwise ($x=1$) set $x=x+1$ and $z=z-1$.
\end{enumerate}
\end{enumerate}
The algorithm is loopless.
Note that the initialization is loopless as well.

The implementation handles all $n\geq$ correctly.
{\codesize
\begin{listing}{1}
// FILE:  src/comb/composition-nz-gray2.h
class composition_nz_gray2
{
    ulong *a_;  // composition: a[1] + a[2] + ... + a[m] = n
    ulong n_;   // compositions of n
    ulong m_;   // current composition has m parts
    ulong e_;   // aux: detection of last composition

    explicit composition_nz_gray2(ulong n)
    {
        n_ = n;
        a_ = new ulong[n_+1+(n_==0)];
        a_[0] = 0;  // returned by last_part() when n==0
        a_[1] = 0;  // returned by first_part() when n==0

        if ( n_ <= 1 )  e_ = n_;
        else            e_ = ( oddq(n_) ? n_ - 1 : n_ );

        first();
    }

    void first()
    {
        if ( n_ <= 1 )
        {
            a_[1] = n_;
            m_ = n_;
        }
        else
        {
            if ( oddq(n_) )
            {
                a_[1] = n_;
                m_ = 1;
            }
            else
            {
                a_[1] = 1;
                a_[2] = n_ - 1;
                m_ = 2;
            }
        }
    }

\end{listing}
}%

A few auxiliary routines are used.
{\codesize
\begin{listing}{1}
    bool oddq(ulong x)  const  { return  0 != ( x & 1UL ); }
    bool evenq(ulong x)  const  { return  0 == ( x & 1UL ); }

    ulong par_to_dir_odd(ulong x)  const
    {
        if ( oddq(x) )  return +1;
        else            return -1UL;
    }

    ulong par_to_dir_even(ulong x)  const
    {
        if ( evenq(x) )  return +1;
        else             return -1UL;
    }
\end{listing}
}%
We split the update into routines
for the cases $z==1$ and $z>1$ as
in the algorithm.
{\codesize
\begin{listing}{1}
    ulong next_zeq1()  // for Z == 1
    {
        const ulong y = a_[m_-1];
        const ulong dy = par_to_dir_odd(y);

        if ( dy == +1 )
        {  // [*, Y, 1 ] --> [*, Y+1 ]
            a_[m_-1] = y + 1;
            m_ -= 1;
            return m_;
        }
        else  // dy == -1
        {  // [*, Y, 1 ] --> [*, Y-1, 1, 1 ]
            a_[m_-1] = y - 1;
            m_ += 1;
            a_[m_] = 1;
            return m_;
        }
    }

    ulong next_zgt1()  // for Z > 1
    {
        const ulong z = a_[m_];
        const ulong y = a_[m_-1];

        if ( oddq(z) )
        {  // [*, Z ] --> [*, Z-1, 1 ]
            a_[m_] = z - 1;
            m_ += 1;
            a_[m_] = 1;
            return m_;
        }

        if ( y != 1 )  // Y > 1
        { // [*, Y, Z ] --> [*, Y+-1, Z-+1 ]
            const ulong dy = par_to_dir_even(y);
            a_[m_-1] = y + dy;
            a_[m_] = z - dy;
            return m_;
        }
        else  // Y == 1
        {
            const ulong x = a_[m_-2];

            if ( x != 1 )
            {  // [*, X, 1, Z ] --> [*, X+-1, 1, Z-+1 ]
                const ulong dx = par_to_dir_odd(x);
                a_[m_-2] = x + dx;
                a_[m_] = z - dx;
                return m_;
            }
            else  // X == 1
            {  // [*, X, 1, Z ] --> [*, X+1, 1, Z-1 ]
                a_[m_-2] = x + 1;
                a_[m_] = z - 1;
                return m_;
            }
        }
    }
\end{listing}
}%
The routine \texttt{next()}
returns number of parts in the new composition
and zero if there are no more compositions.
{\codesize
\begin{listing}{1}
    ulong next()
    {
        ulong z = a_[m_];
        if ( z == e_ )  return 0;  // current is last

        if ( z != 1 )  return  next_zgt1();
        else           return  next_zeq1();
    }
\end{listing}
}%
One update takes about 11 cycles.

We remark that the sublists of compositions into a
fixed number of parts correspond to the combinations
in \jjterm{enup} order.
The list of binary words on the right of figure \ref{fig:comp-gray}
are those shown in \cite[fig.6.6-B (left), p.189]{fxtbook}.

\subsubsection*{Ranking and unranking}

Recursive routines for ranking and unranking
are obtained easily.  The following implementations
are for the ordering with the first part decreasing.
{\codesize
\begin{listing}{1}
// FILE:  src/comb/composition-nz-rank.cc
ulong
composition_nz_gray_rank(const ulong *x, ulong m, ulong n)
// Return rank r of composition x[], 0 <= r < 2**(n-1)
// where n is the sum of all parts.
{
    if ( m <= 1 )  return 0;

    ulong f = x[0];   // first part
    ulong s = n - f;  // remaining sum

    ulong y = composition_nz_gray_rank(x+1, m-1, s);
    if ( 0 == ( f & 1 ) )  // first part even
        return  ( 1UL << (s-1) ) + y;
    else
        return  ( 1UL << s ) - 1 - y;
}

ulong
composition_nz_gray_unrank(ulong r, ulong *x, ulong n)
// Generate composition x[] of n with rank r.
// Return number of parts m of generated composition, 0 <= m <= n.
{
    if ( r == 0 )
    {
        if ( n==0 )  return 0;
        x[0] = n;
        return 1;
    }

    ulong h = highest_one_idx(r);
    ulong f = n - 1 - h;  // first part
    x[0] = f;

    ulong b = 1UL << h;  // highest one
    r ^= b;  // delete highest one

    bool p = f & 1;  // first part f odd ?
    if ( p )  r = b - 1 - r;  // change direction with odd f

    return  1 + composition_nz_gray_unrank( r , x+1, n-f );
}
\end{listing}
}%
%

\section{Appendix I: Nonempty subsets with at most $k$ elements}\label{sect:ksubset}

%
\begin{figure}
{\jjbls
\begin{verbatim}
   1:    1.....    { 0 }                         (continued)
   2:    11....    { 0, 1 }              23:    .1.11.    { 1, 3, 4 }
   3:    111...    { 0, 1, 2 }           24:    .1.1.1    { 1, 3, 5 }
   4:    11.1..    { 0, 1, 3 }           25:    .1..1.    { 1, 4 }
   5:    11..1.    { 0, 1, 4 }           26:    .1..11    { 1, 4, 5 }
   6:    11...1    { 0, 1, 5 }           27:    .1...1    { 1, 5 }
   7:    1.1...    { 0, 2 }              28:    ..1...    { 2 }
   8:    1.11..    { 0, 2, 3 }           29:    ..11..    { 2, 3 }
   9:    1.1.1.    { 0, 2, 4 }           30:    ..111.    { 2, 3, 4 }
  10:    1.1..1    { 0, 2, 5 }           31:    ..11.1    { 2, 3, 5 }
  11:    1..1..    { 0, 3 }              32:    ..1.1.    { 2, 4 }
  12:    1..11.    { 0, 3, 4 }           33:    ..1.11    { 2, 4, 5 }
  13:    1..1.1    { 0, 3, 5 }           34:    ..1..1    { 2, 5 }
  14:    1...1.    { 0, 4 }              35:    ...1..    { 3 }
  15:    1...11    { 0, 4, 5 }           36:    ...11.    { 3, 4 }
  16:    1....1    { 0, 5 }              37:    ...111    { 3, 4, 5 }
  17:    .1....    { 1 }                 38:    ...1.1    { 3, 5 }
  18:    .11...    { 1, 2 }              39:    ....1.    { 4 }
  19:    .111..    { 1, 2, 3 }           40:    ....11    { 4, 5 }
  20:    .11.1.    { 1, 2, 4 }           41:    .....1    { 5 }
  21:    .11..1    { 1, 2, 5 }
  22:    .1.1..    { 1, 3 }
\end{verbatim}
}
\caption{\label{fig:ksubset-lex}
Nonempty subsets of the set $\{0,1,2,\ldots,5\}$ with at most 3 elements.
}
\end{figure}
%

The modifications needed in algorithm \ref{alg:next-sl2}
for restricting the subsets to those with
a prescribed maximal number of elements
are quite small.
Without ado, we give the implementation.
{\codesize
\begin{listing}{1}
// FILE:  src/comb/ksubset-lex.h
class ksubset_lex
{
    ulong n_;  // number of elements in set, should have n>=1
    ulong j_;  // number of elements in subset
    ulong m_;  // max number of elements in subsets
    ulong *x_;  // x[0...j-1]:  subset of {0,1,2,...,n-1}
\end{listing}
}%
The computation of the successor is
{\codesize
\begin{listing}{1}
    ulong next()
    {
        ulong j1 = j_ - 1;
        ulong z1 = x_[j1] + 1;
        if ( z1 < n_ )  // last element is not max
        {
            if ( j_ < m_ )  // append element
            {
                x_[j_] = x_[j1] + 1;
                ++j_;
                return  j_;
            }

            x_[j1] = z1;  // increment last element
            return j_;
        }
        else  // last element is max
        {
            if ( j1 == 0 )  return 0; // current is last

            --j_;
            x_[j_-1] += 1;
            return j_;
        }
    }
\end{listing}
}%
We omit the routine for the predecessor.
Updates take no more time than the updates for the unrestricted subsets.

\section*{Acknowledgment}
I am indebted to Edith Parzefall for
editing two earlier versions of this paper.

\newcommand{\jjbibtitle}[1]{{\small\bfseries #1}}
\newcommand{\jjbibtitletrans}[1]{{\small\bfseries #1}}
\newcommand{\jjbibitem}[2]{\bibitem{#1}{#2}}%
\newcommand{\bdate}[1]{(#1)}



\clearpage
\section{Appendix II: Ranking and unranking methods for mixed radix words}
The following material was added January 2, 2024.

We give C++ code for ranking and unranking mixed radix words
in Subset-lex order and in SL-Gray order.

\subsection*{Subset-lex order}

We put the routines into a C++ class because a table
(variable \texttt{Jmp}) has to be used for ranking and unranking.
It is computed only once, in the constructor.
{\codesize
\begin{listing}{1}
// FILE:  src/comb/mixedradix-subset-lex-rank.h
class mixedradix_subset_lex_rank
// Ranking and unranking function for subset-lex order.
{
    const ulong * A;   // digits of mixed radix number
    const ulong * m1;  // nines for mixed radix base
    const ulong n;     // number of digits
    ulong * Jmp {nullptr};  // jump sizes, see constructor

    mixedradix_subset_lex_rank( const ulong *tA,  // digits
                                ulong tn,         // number of digits
                                const ulong *tm1  // nines
                                )
        : A( tA ), m1( tm1 ), n( tn )
    {
        Jmp = new ulong[n];
        for (ulong j=0; j < n; ++j)
        {
            ulong p = 1;  // product of all radices from index j+1 to end
            for (ulong k=j+1; k < n; ++k)  p *= ( 1 + m1[k] );
            Jmp[j] = p - 1;
        }
    }
\end{listing}
}%

{\codesize
\begin{listing}{1}
    ulong rank()  const
    // Return rank of mixed radix number A[].
    {
        // last index e such that A[e] != 0,
        // e == 0 also for the all-zero word:
        ulong e = 0;
        for (ulong j=0; j<n; ++j)
            if ( A[j] != 0 )  e = j;

        ulong r = 0;
        ulong j = 0;  // track
        while ( true )
        {
            const ulong & aj = A[j];
            const ulong & nine = m1[j];
            const ulong & jmp = Jmp[j];

            if ( j == e )  // on last non-zero digit
            {
                // forward (digits increasing, fast)
                r += aj;
                break;  // and return r
            }
            else  // not on last track
            {
                // forward (digits decreasing, slow)
                r += nine;
                r += jmp * ( nine - aj );
            }

            ++j;
        }

        return r;
    }
\end{listing}
}%

{\codesize
\begin{listing}{1}
    ulong unrank(ulong r, ulong *B)  const
    // Write mixed radix number with rank r into B[].
    // Return index of last non-zero digit, 0 also for the all-zero word.
    {
        ulong j = 0;  // track
        while ( j < n )
        {
            const ulong & nine = m1[j];

            if ( r <= nine )  // on last non-zero digit
            {
                // forward (digits increasing, fast)
                B[j] = r;
                break;
            }
            else  // not on last track
            {
                // forward (digits decreasing, slow)
                r -= nine;
                const ulong & jmp = Jmp[j];
                ulong d = nine;  // new digit

                while ( r > jmp )
                {
                    r -= jmp;
                    d -= 1;
                }

                B[j] = d;
            }

            ++j;
        }

        const ulong tr = j;
        while ( ++j < n )  { B[j] = 0; }
        return tr;
    }
\end{listing}
}%

\subsection*{SL-Gray order}

{\codesize
\begin{listing}{1}
// FILE:  src/comb/mixedradix-sl-gray-rank.h
class mixedradix_sl_gray_rank
// Ranking and unranking function for SL-Gray order.
{
    // Same data and constructor as in mixedradix_subset_lex_rank, omitted.
\end{listing}
}%

{\codesize
\begin{listing}{1}
    ulong rank()  const
    // Return rank of mixed radix number A[].
    {
        // last index e such that A[e] != 0,
        // e == 0 also for the all-zero word:
        ulong e = 0;
        for (ulong j=0; j<n; ++j)
            if ( A[j] != 0 )  e = j;

        ulong r = 0;
        bool fwq = true;  // whether going forward
        ulong j = 0;  // track
        while ( true )
        {
            const ulong & aj = A[j];
            const ulong & nine = m1[j];
            const ulong & jmp = Jmp[j];

            if ( j == e )  // on last non-zero digit
            {
                if ( fwq )  // forward (digits increasing, fast)
                {
                    r += aj;
                }
                else        // backward (digits decreasing, fast)
                {
                    const ulong ten = nine + 1;
                    r += jmp * ten;
                    r += ( ten - aj );
                }

                break;  // and return r
            }
            else  // not on last track
            {
                if ( fwq )  // forward (digits decreasing, slow)
                {
                    r += nine;
                    r += jmp * ( nine - aj );
                }
                else        // backward (digits increasing, slow)
                {
                    r += jmp * aj;
                }
            }

            // may switch direction:
            const ulong ap = aj & 1UL;
            const ulong mp = nine & 1UL;
            if ( ap != mp )  fwq = ! fwq;

            ++j;
        }

        return r;
    }
\end{listing}
}%

{\codesize
\begin{listing}{1}
    ulong unrank(ulong r, ulong *B)  const
    // Write mixed radix number with rank r into B[].
    // Return index of last non-zero digit, 0 also for the all-zero word.
    {
        bool fwq = true;  // whether going forward

        ulong j = 0;  // track
        while ( j < n )
        {
            const ulong & nine = m1[j];
            const ulong & jmp = Jmp[j];

            if ( fwq )  // forward (same as for subset-lex order)
            {
                if ( r <= nine )   // forward (digits increasing, fast)
                {
                    // on last non-zero digit
                    B[j] = r;
                    break;  // done
                }
                else               // forward (digits decreasing, slow)
                {
                    r -= nine;
                    ulong d = nine;  // new digit
                    while ( r > jmp )
                    {
                        r -= jmp;
                        d -= 1;
                    }
                    B[j] = d;
                }
            }
            else  // backward
            {
                const ulong ten = nine + 1;
                const ulong jj = jmp * ten;
                if ( r > jj )      // backward (digits decreasing, fast)
                {
                    // on last non-zero digit
                    r -= jj;
                    B[j] = ten - r;
                    break;  // done
                }
                else               // backward (digits increasing, slow)
                {
                    ulong d = 0;  // new digit
                    while ( r > jmp )
                    {
                        r -= jmp;
                        d += 1;
                    }
                    B[j] = d;
                }
            }

            const ulong ap = B[j] & 1UL;
            const ulong mp = nine & 1UL;
            if ( ap != mp )  fwq = ! fwq;  // switch direction

            ++j;
        }

        const ulong tr = j;
        while ( ++j < n )  { B[j] = 0; }
        return tr;
    }
\end{listing}
}%

\end{document}